\documentclass [letterpaper,11pt]{article}

\usepackage[T1]{fontenc}
\usepackage[utf8]{inputenc}
\usepackage[english]{babel}
\usepackage[margin=1in]{geometry}

\usepackage{amsmath}
\usepackage{amssymb}
\usepackage{amsfonts}
\usepackage{array}
\usepackage{float}
\usepackage{graphicx}
\usepackage[font={footnotesize}]{caption}
\usepackage{subcaption}
\usepackage{appendix}
\usepackage{xfrac}
\usepackage{placeins}
\usepackage{bm}
\usepackage{cancel}
\usepackage{amsthm}
\usepackage[dvipsnames]{xcolor}
\usepackage{comment}
\usepackage{booktabs} 
\usepackage{physics}
\usepackage{mathtools}

\usepackage[hidelinks]{hyperref}
\hypersetup{colorlinks=true, citecolor=blue, linkcolor=blue, urlcolor=blue}

  \definecolor{ICES}{RGB}{94, 156, 174}
  \definecolor{ORANGE}{RGB}{191, 87, 0}
  \definecolor{RED}{RGB}{190, 30, 49}
  \definecolor{SUN}{RGB}{227, 81, 51}
  \definecolor{GREEN}{RGB}{0, 171, 86}
  \definecolor{BLUE}{RGB}{11, 78, 179}
  \definecolor{BROWN}{RGB}{122, 80, 40}
  \definecolor{GREY}{RGB}{50, 50, 50}
  \definecolor{TEAL}{RGB}{0, 160, 176}

\newtheorem*{remark}{\emph{Remark}}

\numberwithin{equation}{section}

\newcommand{\be}{\begin{equation}}
\newcommand{\ee}{\end{equation}}

\newcommand{\eq}[1]{(\ref{eq:#1})}
\newcommand{\fig}[1]{\ref{fig:#1}}
\newcommand{\tab}[1]{\ref{tab:#1}}

\newcommand{\bb}[1]{\mathbb{#1}}
\newcommand{\mc}[1]{\mathcal{#1}}

\newcommand{\p}{\partial}

\newcommand{\bfinvnoise}{\mathbf{\Gamma}_{\hskip -2pt \text{noise}}^{-1}}

\newcommand{\noisevec}{\boldsymbol{\nu}}

\newcommand{\bfnoise}{\mathbf{\Gamma}_{\hskip -2pt \text{noise}}}

\newcommand{\bfpostcov}{\mathbf{\Gamma}_{\hskip -2pt \text{post}}}

\newcommand{\piprior}{\pi_{\text{prior}}}
\newcommand{\pipost}{\pi_{\text{post}}}
\newcommand{\pilike}{\pi_{\text{like}}}

\newcommand{\bfmmap}{\vb{m}_{\text{map}}}
\newcommand{\bfmprior}{\vb{m}_{\text{prior}}}
\newcommand{\bfmtrue}{\vb{m}_{\text{true}}}
\newcommand{\bfmmapprior}{\vb{m}_{\text{map}}^{\text{prior}}}

\newcommand{\diag}{\text{diag}}

\newcommand{\bfpostcovm}{\mathbf{\Gamma}_{\hskip -2pt \text{post}(m)}}
\newcommand{\bfpriorcovm}{\mathbf{\Gamma}_{\hskip -2pt \text{prior}(m)}}

\newcommand{\bfinvpriorcovm}{\mathbf{\Gamma}_{\hskip -2pt \text{prior}(m)}^{-1}}
\newcommand{\bfinvpostcovm}{\mathbf{\Gamma}_{\hskip -2pt \text{post}(m)}^{-1}}

\newcommand{\bfpostcovq}{\mathbf{\Gamma}_{\hskip -2pt \text{post}(q)}}
\newcommand{\bfpriorcovq}{\mathbf{\Gamma}_{\hskip -2pt \text{prior}(q)}}

\newcommand{\bfqmap}{\vb{q}_{\text{map}}}
\newcommand{\bfqprior}{\vb{q}_{\text{prior}}}
\newcommand{\bfqtrue}{\vb{q}_{\text{true}}}
\newcommand{\bfqmapprior}{\vb{q}_{\text{map}}^{\text{prior}}}

\def\hH1{H^1(\mathcal{T}_h)}

\def\H1{H^1(\Omega)}

\DeclareMathOperator*{\argmin}{arg\,min}

\newcommand{\paramtxt}{m}
\newcommand{\datatxt}{d}
\newcommand{\paramvec}{\vb{\paramtxt}}
\newcommand{\datavec}{\vb{\datatxt}}
\newcommand{\obsvec}{\vb{\datatxt}^{\text{obs}}}
\newcommand{\statevec}{\vb{w}}
\newcommand{\qoivec}{\vb{q}}

\newcommand{\Imat}{\mathbf{I}}
\newcommand{\Kmat}{\mathbf{K}}
\newcommand{\Qmat}{\mathbf{Q}}
\newcommand{\ptomat}{\mathbf{F}}
\newcommand{\ptqmat}{\mathbf{F}_{\hskip -1pt q}}
\newcommand{\Gmat}{\mathbf{G}}

\newcommand{\numstate}{N_w}
\newcommand{\numtime}{N_t}
\newcommand{\numparam}{N_m}
\newcommand{\numdata}{N_d}
\newcommand{\numqoi}{N_q}

\usepackage[backend=bibtex,
	maxbibnames=99,
	maxcitenames=99,
	giveninits=true,
	doi=true,
	isbn=false]
	{biblatex}
\graphicspath{{./}}

\title{\Large{
Goal-Oriented Real-Time Bayesian Inference for\\
Linear Autonomous Dynamical Systems\\
With Application to Digital Twins for Tsunami Early Warning
}}
\author{
Stefan Henneking$^{1,*}$,
Sreeram Venkat$^{1}$,
Omar Ghattas$^{1,2}$ \\[2pt]
{\small $^{1}$Oden Institute for Computational Engineering and Sciences,
The University of Texas at Austin} \\
{\small $^{2}$Walker Department of Mechanical Engineering, 
The University of Texas at Austin}
}
\date{\today}

\begin{document}

\clearpage \maketitle
\thispagestyle{empty}
\renewcommand*{\thefootnote}{\fnsymbol{footnote}}
\footnotetext{$^*$Corresponding author: stefan@oden.utexas.edu}
\renewcommand*{\thefootnote}{\arabic{footnote}}

\addcontentsline{toc}{section}{Abstract}

\begin{abstract}
  We present a goal-oriented framework for constructing digital twins
  with the following properties: (1)~they employ discretizations of
  high-fidelity partial differential equation (PDE) models governed by
  autonomous dynamical systems, leading to large-scale forward
  problems; (2)~they solve a linear inverse problem to assimilate
  observational data to infer uncertain model components followed by a
  forward prediction of the evolving dynamics; and (3)~the entire
  end-to-end, data-to-inference-to-prediction computation is carried
  out without approximation and in real time through a Bayesian
  framework that rigorously accounts for uncertainties. Several
  challenges must be overcome to realize this framework, including the
  large scale of the forward problem, the high dimensionality of the
  parameter space, and for a class of problems including those we target,
  the slow decay of the
  singular values of the parameter-to-observable map. Here we
  introduce a methodology to overcome these challenges by exploiting
  the autonomous structure of the forward model to decompose the
  solution of the inverse problem into a one-time-only offline phase
  in which the PDE model is solved a limited number of times (equal to
  the number of sensors), and an online phase that maps well onto GPUs
  and computes the parameter inference and prediction of quantities of
  interest in real time, given observational data.
  Our ultimate goal is to apply this framework to construct digital
  twins for subduction zones, including Cascadia, to provide early
  warning for tsunamis generated by
  megathrust earthquakes. To this end, we demonstrate how our
  methodology can be used to employ seafloor pressure observations,
  along with the coupled acoustic--gravity wave equations, to infer
  the earthquake-induced spatiotemporal seafloor motion (discretized
  with $\mathcal{O}(10^9)$ parameters) and forward predict the tsunami
  propagation. We present results of an end-to-end inference,
  prediction, and uncertainty quantification for a representative test
  problem with $\mathcal{O}(10^8)$ inversion parameters for which 
  goal-oriented Bayesian inference is accomplished exactly and in real
  time, that is, in a matter of seconds.
\end{abstract}

\vspace{10pt}


\section{Introduction}
\label{sec:intro}

\paragraph{Motivation}

This paper presents a novel framework for real-time solution of very
large scale linear inverse problems governed by autonomous dynamical
systems governed by high-fidelity partial differential equation (PDE)
models, with quantified uncertainties via the framework of Bayesian
inference. Our target application is the inference of
megathrust earthquake-induced seafloor motion from seafloor
measurements in subduction zones, and subsequent forecasting of
tsunamis in vulnerable coastal regions, with particular focus on the
Cascadia subduction zone (CSZ) in the Pacific Northwest.  Tsunamis
resulting from megathrust earthquakes are capable of massive
destruction and loss of human life.  Two recent examples are the
tsunamis unleashed by the 2011 T\=ohoku earthquake, resulting in
22,000 deaths and $\$410$~B in damages, and the 2004 Sumatra--Andaman
earthquake, resulting in 230,000--280,000 total deaths.  Strain has
been accumulating on the CSZ since its last major rupture in 1700.
Statistical recurrence intervals suggest a significant probability of
a major rupture (up to magnitude 9) on this fault in the next 50
years~\cite{leonard2010rupture}, which would have an enormous impact
on millions of people living in the Pacific
Northwest~\cite{cherniawsky2007tsunami, olsen2008earthquake}.

With the increasingly dense current or planned instrumentation of
subduction zones, including Cascadia~\cite{schmidt2019monitoring,
  becerril2024towards}, and the advance of earthquake--tsunami
simulations enabled by high performance computing~\cite{uphoff2017sumatra, krenz2021palu, abrahams2023comparison}, the availability of real-time observational data promises to enhance tsunami early warning systems.
Our methodology aims to employ pressure observations obtained from
seafloor sensors, along with the coupled acoustic--gravity wave
equations that describe the generation of
tsunamis~\cite{yamamoto1982gravity, lotto2015tsunami}, to infer the
earthquake-induced spatiotemporal seafloor motion and quantify the
uncertainty in its inference. 
Employing a model that accounts for acoustic wave propagation is critical for extracting the maximum amount of information from the acoustic pressure sensors, particularly in the near field and during the transients that follow a rupture on the fault~\cite{leveque2018cascadia, li2009ocean}.
The solution of this inverse problem then provides the seafloor forcing to forward propagate the tsunamis toward populated and vulnerable regions along coastlines.
Other tsunami early warning methods that quantify uncertainties in tsunamigenic earthquake sources (e.g.~\cite{giles2021faster} and references therein) or directly make predictions of tsunami waveforms (e.g.~\cite{fujita2024scenario} and references therein) typically rely on statistical tools (rather than full-order PDE solutions) to enable rapid forecasting of the propagating tsunami under uncertainty.
To provide sufficient early warning to guide evacuation of lower lying regions and inform early responders, the entire inverse problem and subsequent forecasting of coastal inundation must be carried out in real time, i.e.~in a few seconds.
This is because tsunamis generated on the CSZ may reach coastal communities in as little as 20 minutes~\cite{leveque2018cascadia, melgar2016kinematic}.
Such early warning systems would have saved innumerable lives in the T\=ohoku and Sumatra--Andaman earthquakes.
However, solution of this inverse problem, governed by the discretized
high-fidelity PDE model, with quantified uncertainties and in real time,
appears to be intractable. 

\paragraph{Large-scale inverse problems}

An appropriate discretization of the acoustic--gravity PDE forward
model for the $1000$~km long CSZ gives rise to a system with
$\mathcal{O}(10^{9})$ inversion parameters describing the
spatiotemporal earthquake-induced seafloor motion.  Scalable methods
for large-scale Bayesian inverse problems governed by PDEs with
heterogeneous 3D inversion parameter fields
have been
developed and successfully employed for a number of challenging
applications, including 
ice sheet flow~\cite{isaac2015scalable},
global seismology~\cite{Bui-ThanhGhattasMartinEtAl13},
poroelasticity~\cite{AlghamdiHesseChenEtAl21},
atmospheric transport~\cite{flath2011fast}, and
ocean dynamics~\cite{kalmikov2014ocean}. 
Typically, these methods rely on exploiting the low-dimensional
structure of the map from parameters (inputs) to observations
(outputs) in order to facilitate scalable sampling from the posterior,
or a Laplace approximation of it. 
Alternatively, surrogates of the map from parameters to either the observables or
the data misfit functional can be constructed and used as proxies
for the forward problem. Various surrogate methods exist, 
including Gaussian
processes~\cite{kennedy2001bayesian}, radial basis
functions~\cite{wild2008orbit}, projection-based model
reduction~\cite{benner2015survey},
and polynomial chaos expansions
\cite{MarzoukNajmRahn07,MarzoukXiu09},
as well as neural network-based surrogates such as Fourier neural
operators~\cite{li2021fourier}, DeepONets~\cite{LuJinPangEtAl21},
and derivative-informed neural
operators~\cite{oleary2024dino,CaoOLearyRoseberryGhattas2025}. 
Generally speaking, these techniques work well when the
parameter-to-observable map can be well-approximated in a
lower-dimensional
manifold, which is indeed a feature inherent in many problems in
science and engineering~\cite{ghattas2021learning}.

Unfortunately, the tsunami inverse problem does not admit
a low-dimensional subspace representation~\cite{greif2019kolmogorov}
due to the governing high-frequency wave propagation forward problem. 
Moreover, wave propagation inverse problems are ill-posed in the high
wavenumber components of the parameter field \cite{bui2012a, bui2012b,
  bui2014} due to band-limited observational data; additional
``noise'' stemming from the surrogate (due to its typical inability to
capture high frequencies) will be amplified in the inverse solution
and may lead to instability, even if instrument noise is low. 
The high-rank structure of the input--output map and the corresponding
inverse operator---i.e.~the Hessian matrix of the negative log
posterior---renders existing techniques
inefficient and therefore inapplicable to real-time inference of
high-dimensional parameter fields for these types of inverse problems. 
The need to forward predict specific quantities of interest
(QoIs) and quantify their uncertainties within a goal-oriented
Bayesian framework further compounds the difficulties of the problem. 
In recent years, several methods for efficiently approximating high-rank Hessians have
been developed, including those that exploit the
pseudo-differential~\cite{nammour2011phd, demanet2012probing},
augmented Lagrangian~\cite{alger2017augmented},
product-convolution~\cite{alger2019product},
H-matrix~\cite{ambartsumyan2020hierarchical}, and point spread
function~\cite{alger2024point} structure of the Hessians of particular
classes of inverse problems. Despite these advances, the $10^9$-dimensional problems
we target remain out of reach. 


Rather, we can entirely circumvent the need to approximate the
Hessian, and instead exactly solve the (discretized) inverse problem
when it is governed by an autonomous dynamical system.  The evolution
of such systems with respect to any given input may depend on the
system’s current state but does not explicitly depend on the
independent variable.  The coupled acoustic--gravity PDE model that
governs the tsunami dynamics is a time-invariant dynamical system,
which is a subclass of autonomous systems where the independent
variable is time.  For such problems, we can exploit this problem
structure and move most of the computationally expensive parts of the
Bayesian inverse problem offline.  As a key component of this
framework, we employ a multi-GPU-accelerated fast Fourier transform
(FFT)-based algorithm for computing Hessian matrix--vector products
(matvecs) that provides a $\mc{O}(10^4)\times$ speedup over
conventional methods (after a one-time setup
cost)~\cite{venkat2024fft}.

\paragraph{Contributions}

Solution of inverse problems
in actionable time scales followed by prediction of QoIs under
uncertainty is a fundamental component of digital twins \cite{Asch22,
  RasheedSanKvamsdal20, DigitalTwinsNASEM,
  NiedererSacksGirolamiEtAl21}.  This paper introduces an efficient
and scalable digital twin framework to perform goal-oriented real-time
Bayesian inference for linear autonomous dynamical systems. The
framework is applicable to inverse problems with high-rank Hessians
for which existing surrogate techniques are inefficient.  Moreover,
contrary to methods based on surrogate models, our approach solves the
inverse problem exactly, i.e.~the inverse solution exactly satisfies
the (discretized) PDE model and the maximum-a-posterior point is
computed exactly.  The framework can be used to infer
infinite-dimensional parameter fields (high-dimensional after
discretization) but assumes that the observational data are relatively
sparse. In particular, the number of spatial observers (e.g.~sensors)
is assumed to be small enough such that solution of the governing PDEs
once per sensor location is feasible.  The Bayesian inverse problem is
then posed in a goal-oriented setting that forward predicts
QoIs---including their uncertainties---once the parameter field is
inferred.

Our digital twin framework for goal-oriented real-time Bayesian
inference is made possible by the following key ideas: (1)~an
offline--online decomposition of the inverse solution that precomputes
various mappings to enable the real-time application of the inverse
operator and entirely circumvents the need for PDE solutions at the
time of inference; (2)~representation of the inverse operator in the
lower-dimensional data space instead of the high-dimensional parameter space;
(3)~extension of this real-time inference methodology to the
goal-oriented setting for the prediction of QoIs under uncertainty;
(4)~exploitation of the shift invariance of autonomous dynamical
systems to reduce the required number of PDE solutions by several
orders of magnitude; and (5)~exploitation of this same property to
perform efficient FFT-based and multi-GPU-accelerated Hessian matvecs.
We show how this framework can be applied to construct digital twins
for tsunami early warning from sea bottom pressure observations and a
coupled acoustic--gravity PDE model.  For a representative test
problem with more than $10^8$ parameters and 800~QoIs, both the
parameter inference and the full end-to-end,
data-to-inference-to-prediction computations are carried out within a
fraction of a second.

\section{Background}
\label{sec:background}

\paragraph{Notation}
Lower-case italic letters ($a$) are used to denote infinite-dimensional functions or variables.
Discretized functions (vectors in $\bb{R}^n$) are denoted by lower-case upright boldface ($\vb a$) notation.
Upper-case calligraphic letters ($\mc A$) indicate infinite-dimensional operators.
Discretized operators (matrices in $\bb{R}^{m \times n}$) are indicated by upper-case upright boldface ($\vb A$) notation.
Lower-case italic ($c$) or upper-case italic ($C$) letters may also denote constants.
Greek letters are used in different contexts, but their definitions will be made clear as they appear.
Any exceptions to the above will be clear from context.
The notation $\| \vb a \|_{\vb A}$ refers to the weighted Euclidean inner product $\vb a^T \! \vb A \vb a$; if no weight is given, the weight is identity.
Estimates for computation times are given in h~(hours), m~(minutes), s~(seconds), or ms~(milliseconds) of wall-clock time.
Storage estimates are based on double precision (i.e.~eight bytes per floating point number) and given in PB~(Petabytes), TB~(Terabytes), GB~(Gigabytes), or MB~(Megabytes).

\paragraph{Linear time-invariant dynamical system}
While our framework applies to linear autonomous dynamical systems in general, we consider the case of a linear time-invariant (LTI) dynamical system, which is applicable to the tsunami early warning application discussed in this paper. In particular, we consider an LTI dynamical system of the form
\be
\begin{aligned}
   \frac{\partial w}{\partial t} &= \mc A w + \mc C m && \text{in } \Omega \times (0,T), \\
   w &= w_0 && \text{in } \Omega \times \{ 0 \}, \\
   d &= \mc B w && \text{in } \Omega \times (0,T),
\end{aligned}
\label{eq:LTI}
\ee
with appropriate boundary conditions on $\p \Omega \times (0,T)$, where $\Omega$ is the spatial domain, $(0,T)$ is the time domain, $w(x,t)$ is the system's state with initial state $w_0(x)$, parameter (input) $m(x,t)$ represents the source or forcing of the system and is independent of the state, and both $\mc A$ and $\mc C$ are time-invariant differential operators; $d(x,t)$ describes the observables (output) of the system, which are extracted from the state $w$ via a time-invariant observation operator $\mc B$.

\paragraph{Discrete LTI system}
Consider a discrete version of the LTI system obtained by discretizing \eq{LTI} in time with a single-step explicit method\footnote{Note that the methodology is easily extended to multi-step explicit and implicit methods.} and with some spatial discretization scheme,
\be
   \statevec_{k+1} = \vb A \statevec_k + \vb C \paramvec_k, \quad k=0,1,\cdots,\numtime-1,
   \label{eq:LTI-time-stepping}
\ee
where $\statevec_k \in \bb{R}^{\numstate}$, $\paramvec_k \in \bb{R}^{\numparam}$, and the discrete operators $\vb A \in \bb{R}^{\numstate \times \numstate}$ and $\vb C \in \bb{R}^{\numstate \times \numparam}$ both depend on the particular time-stepping scheme and spatial discretization method. 
Then, using \eq{LTI-time-stepping} the discretized LTI system can be written as follows:

\be
\begin{split}
   \statevec_1 &= \vb A \statevec_0 + \vb C \paramvec_0 , \\
   \statevec_2 &= \vb A \statevec_1 + \vb C \paramvec_1 , \\
      &= \vb A (\vb A \statevec_0 + \vb C \paramvec_0) + \vb C \paramvec_1 = \vb A^2 \statevec_0 + \vb A^1 \vb C \paramvec_0 + \vb A^0 \vb C \paramvec_1 , \\[-2pt]
      \vdots \\[-10pt]
   \statevec_{k+1} &= \vb A^{k+1} \statevec_0 + \sum_{i=0}^k \vb A^i \vb C \paramvec_{k-i} , \\[-3pt]
   \datavec_{k+1} &= \vb B \statevec_{k+1} ,
\end{split}
\label{eq:LTI-time-discrete}
\ee
where $\datavec_k \in \bb{R}^{\numdata}$ and $\vb B \in \bb{R}^{\numdata \times \numstate}$ is the discrete observation operator.

Without loss of generality, assume homogeneous initial condition $\vb u_0 = \vb 0$. We can then write the discretized LTI system in the following way:

\be
   \left[ 
   \begin{array}{@{\hskip 2pt}c@{\hskip 2pt}}
      \datavec_1 \\
      \datavec_2 \\
      \vdots \\
      \datavec_{k+1} \\
      \vdots \\
      \datavec_{\numtime}
   \end{array}
   \right]
   =
   \left[ 
   \begin{array}{@{\hskip 2pt}cc@{\hskip 4pt}c@{\hskip 4pt}c@{\hskip 4pt}c@{\hskip 4pt}c@{\hskip 2pt}}
      \vb B \vb A^0 \vb C  \\
      \vb B \vb A^1 \vb C & \vb B \vb A^0 \vb C  \\
      \vdots & \vdots & \ddots  \\
      \vb B \vb A^k \vb C & \vb B \vb A^{k-1} \vb C & \cdots & \vb B \vb A^0 \vb C \\
      \vdots & \vdots & & \vdots & \ddots  \\
      \vb B \vb A^{\numtime-1} \vb C & \vb B \vb A^{\numtime-2} \vb C & \cdots & \vb B \vb A^{\numtime-(k+1)} \vb C & \cdots & \vb B \vb A^0 \vb C
   \end{array}
   \right]
   \left[ 
   \begin{array}{@{\hskip 2pt}c@{\hskip 2pt}}
      \paramvec_0 \\
      \paramvec_1 \\
      \vdots \\
      \paramvec_k \\
      \vdots \\
      \paramvec_{\numtime-1}
   \end{array}
   \right] .
\label{eq:LTI-matrix}
\ee

We define $\ptomat_{ij} \coloneqq \vb B \vb A^{i-j} \vb C \in \bb{R}^{\numdata \times \numparam}$, $i,j=1,2,\ldots,\numtime, i \ge j$; the LTI system \eq{LTI-time-stepping} can then be written more compactly as follows:

\begin{equation}
   \left[ \begin{array}{c}
   \datavec_1 \\[2pt]
   \datavec_2 \\[6pt]
   \datavec_3 \\[1pt]
   \vdots \\[3pt]
   \datavec_{\numtime}
   \end{array} \right]
   =
   \left[ \begin{array}{ccccc}
   \ptomat_{11} & \vb 0 & \vb 0 & \cdots & \vb 0 \\[2pt]
   \ptomat_{21} & \ptomat_{11} & \vb 0 & \cdots & \vb 0 \\
   \ptomat_{31} & \ptomat_{21} & \ptomat_{11} & \ddots & \vdots \\
   \vdots & \vdots & \ddots & \ddots & \vb 0 \\[2pt]
   \ptomat_{\numtime,1} & \ptomat_{\numtime-1,1} & \cdots & \ptomat_{21} & \ptomat_{11}
   \end{array} \right]
   \hskip 5pt
   \left[ \begin{array}{c}
   \paramvec_0 \\[2pt]
   \paramvec_1 \\[6pt]
   \paramvec_2 \\[1pt]
   \vdots \\[3pt]
   \paramvec_{\numtime-1}
   \end{array} \right] ,
   \label{eq:ShiftInvariance}
\end{equation}
or very concisely as
\be
   \datavec \coloneqq \ptomat \paramvec .
   \label{eq:p2o}
\ee

We refer to $\ptomat$ as the (discrete) \emph{parameter-to-observable} (p2o) map; $\paramvec$ is the parameter vector and $\datavec$ is the vector of observables or data vector.
Then,

\begin{itemize}
   \item $\paramvec \in \bb R^{\numparam \numtime}$ with blocks $\paramvec_j \in \bb R^{\numparam}$, $j = 0,1, \ldots, \numtime-1$;
   \item $\datavec \in \bb R^{\numdata \numtime}$ with blocks $\datavec_i \in \bb R^{\numdata}$, $i = 1,2, \ldots, \numtime$;
   \item $\ptomat \in \bb R^{(\numdata \numtime) \times (\numparam \numtime)}$ with blocks $\ptomat_{ij} \in \bb{R}^{\numdata \times \numparam}$, $i,j = 1,2, \ldots, \numtime$.
\end{itemize}

It is clear from \eq{ShiftInvariance} that the p2o map $\ptomat$ is shift-invariant with respect to its blocks $\ptomat_{ij}$. In particular, $\ptomat$ is \emph{block Toeplitz}. Additionally, time causality implies that $\ptomat$ is block lower-triangular. Section~\ref{sec:methodology} describes how this special structure of $\ptomat$ can be effectively exploited in the context of solving inverse problems.

\paragraph{Bayesian inverse problem}
Given a system of the form \eq{p2o}, we consider the inverse problem of inferring the parameters $\paramvec$ from the observed data $\obsvec$.
In finite dimensions, Bayes' theorem is given by:
\be
   \pipost(\paramvec | \obsvec) \propto \pilike(\obsvec | \paramvec) \piprior(\paramvec) ,
\ee
where $\piprior(\paramvec)$ is the prior probability density of the model parameters, $\pilike(\obsvec | \paramvec)$ is the likelihood of the data $\obsvec$ given parameters $\paramvec$, and $\pipost(\paramvec | \obsvec)$ denotes the posterior density reflecting the probability of the parameters conditioned on the data.

We assume a Gaussian prior
\be
   \piprior(\paramvec) \propto
   \exp \left\{
      - \frac{1}{2} \| \paramvec - \bfmprior \|_{\bfpriorcovm^{-1}}^2
   \right\} ,
\ee
where $\bfmprior$ and $\bfpriorcovm$ are respectively the mean and covariance of the prior distribution.

Given data $\obsvec$ perturbed by additive Gaussian noise, i.e.~$\obsvec = \ptomat \paramvec + \noisevec$, $\noisevec \sim \mc N(\vb 0, \bfnoise)$, the likelihood is
\be
   \pilike(\obsvec | \paramvec) \propto
   \exp \left\{ 
   - \frac{1}{2} \| \ptomat \paramvec - \obsvec \|_{\bfnoise^{-1}}^2
   \right\} .
\ee

For this choice of Gaussian prior and additive Gaussian noise, the posterior is then given by:
\be
   \pipost(\paramvec  | \obsvec) \propto
   \exp \left\{
      - \frac{1}{2} \| \ptomat \paramvec - \obsvec \|_{\bfnoise^{-1}}^2
      - \frac{1}{2} \| \paramvec - \bfmprior \|_{\bfpriorcovm^{-1}}^2
   \right\} .
\ee

The maximum-a-posterior (MAP) point $\bfmmap$ maximizes the posterior distribution over admissible parameters $\paramvec$:
\be
   \bfmmap \coloneqq \argmin_{\paramvec \in \bb{R}^{\numparam \numtime}} (-\log \pipost(\paramvec))
   = \argmin_{\paramvec \in \bb{R}^{\numparam \numtime}} 
   \frac{1}{2} \| \ptomat \paramvec - \obsvec \|_{\bfnoise^{-1}}^2 +
   \frac{1}{2} \| \paramvec - \bfmprior \|_{\bfpriorcovm^{-1}}^2 .
\ee

The MAP point $\bfmmap$ is thus given by the solution to the linear inverse problem
\be
   \left( \ptomat^* \bfnoise^{-1} \ptomat + \bfpriorcovm^{-1} \right) \bfmmap =
   \ptomat^* \bfnoise^{-1} \obsvec + \bfpriorcovm^{-1} \bfmprior,
   \label{eq:mmap-problem}
\ee
where $\ptomat^*$ denotes the adjoint of the p2o map and $\vb H \coloneqq \left( \ptomat^* \bfnoise^{-1} \ptomat + \bfpriorcovm^{-1} \right)$ is the Hessian of the negative log-posterior evaluated at $\bfmmap$.

The posterior for this linear inverse problem is then a Gaussian centered at the MAP point, i.e.~$\pipost(\paramvec | \obsvec) \propto \mc N(\bfmmap, \bfpostcovm)$, with posterior covariance
\be
   \bfpostcovm \coloneqq \vb H^{-1} = \left( \ptomat^* \bfnoise^{-1} \ptomat + \bfpriorcovm^{-1} \right)^{-1} .
   \label{eq:Hessian}
\ee

For a more thorough discussion of Bayesian inverse problems, see \cite{ghattas2021learning, villa2021hippylib} and references therein.

\paragraph{Goal-oriented Bayesian inference}
In the goal-oriented setting, we consider an extension to the LTI system \eq{LTI}:
\be
   q = \mc B_q w \text{ in } \Omega \times (0,T) ,
   \label{eq:LTI-qoi}
\ee
where some QoI, denoted by $q$, is extracted from the system's state via a time-invariant QoI observation operator $\mc B_q$.

Analogous to \eq{p2o}, the discrete LTI system for the QoI can be defined concisely by
\be
   \qoivec \coloneqq \ptqmat \paramvec ,
   \label{eq:p2q}
\ee
where $\ptqmat$ is the discrete \emph{parameter-to-QoI} (p2q) map and $\qoivec$ is the QoI vector; that is,
\begin{itemize}
   \item $\qoivec \in \bb R^{\numqoi \numtime}$ with blocks $\qoivec_i \in \bb R^{\numqoi}$, $i = 1,2, \ldots, \numtime$;
   \item $\ptqmat \in \bb R^{(\numqoi \numtime) \times (\numparam \numtime)}$ with blocks $(\ptqmat)_{ij} \in \bb{R}^{\numqoi \times \numparam}$, $i,j = 1,2, \ldots, \numtime$.
\end{itemize}
Similar to the block-wise definition of the p2o map $\ptomat$ in \eq{ShiftInvariance}, $\ptqmat$ is defined by its
blocks, 
\be
   (\ptqmat)_{ij} \coloneqq \vb{B}_q \vb A^{i-j} \vb C \in \bb{R}^{\numqoi \times \numparam},\quad i,j=1,2,\ldots,\numtime, i \ge j,
\ee
where $\vb{B}_q \in \bb{R}^{\numqoi \times \numstate}$ is the discrete QoI observation operator.
Like the p2o map $\ptomat$, the p2q map $\ptqmat$ is block lower-triangular and shift-invariant with respect to its block, 
i.e.~$\ptqmat$ is also block Toeplitz.

\begin{remark}
In general, the number of time steps used by the time-stepping algorithm ($\numtime$) may be different from the temporal discretization of the parameters $(\numtime^m)$, data $(\numtime^d)$, and QoIs $(\numtime^q)$.
The block structure of the p2o and p2q maps then changes correspondingly but the same algorithmic techniques apply.
For ease of presentation, we use $\numtime = \numtime^m = \numtime^d = \numtime^q$ here.
The numerical example in Section~\ref{sec:numerical} presents a case where $\numtime^q \ll \numtime^d = \numtime^m \ll \numtime$.
\end{remark}

Now, consider the problem of predicting the QoIs $\qoivec$ from the observed data $\obsvec$ via the inference of parameters $\paramvec$.
Since the Bayesian prior and posterior of the parameter field are Gaussians, and a linear transformation of a Gaussian is still Gaussian, we can explicitly write down the prior and posterior for $\qoivec$:
\begin{align}
   \piprior(\qoivec) &\propto \mc N \left( \bfqprior, \bfpriorcovq \right) 
   = \mc N \left( \ptqmat \bfmprior, \ptqmat \bfpriorcovm \ptqmat^* \right) , \\[3pt]
   \pipost(\qoivec | \obsvec) &\propto \mc N \left( \bfqmap, \bfpostcovq \right)
   = \mc N \left( \ptqmat \bfmmap, \ptqmat \bfpostcovm \ptqmat^* \right) ,
   \label{eq:qoi-posterior}
\end{align}
where $\ptqmat^*$ is the adjoint p2q map.
In other words, the means of the prior and posterior of the QoIs can be computed from the means of the prior and posterior of the inferred parameters via a push-forward through the p2q map; the respective covariances require an additional application of the adjoint p2q map.

Ultimately, goal-oriented Bayesian inference aims to estimate predictions of the QoIs from observed data without the need for inference of the parameters. In Section~\ref{sec:methodology}, we discuss how this can be accomplished in our mathematical and computational framework for real-time inference and prediction.

\section{Tsunami Inverse Model}
\label{sec:model}

The Bayesian framework for real-time inference and prediction
presented in this paper is applicable to a variety of application
problems that are governed by PDE operators with underlying autonomous
structure.  Some examples will be given in the concluding remarks.
However, rather than introducing this framework abstractly, we believe
it is more descriptive to the reader to follow the presentation along
with a concrete application problem in mind.  For this reason, we
introduce our motivating application in this section before discussing
the framework step-by-step in the subsequent section, where this
application will serve as an example both in the description of the
framework and in giving concrete estimates of the memory and
computational complexities involved in realizing it for a large-scale
application (even if those estimates may differ considerably for other
application problems).

\paragraph{Acoustic--gravity wave equations}
The tsunami dynamics are modeled by the linear acoustic--gravity wave equations that couple the propagation of ocean acoustic waves with the surface gravity wave.
This PDE model is derived by linearization of the conservation of mass and momentum around hydrostatic pressure in the compressible ocean, and the coupling to the surface gravity wave comes from a modified free surface boundary condition at the sea surface~\cite{lotto2015tsunami}.
The model assumes that the surface gravity wave height is much smaller
than the ocean depth, so the model may be a good approximation in
the ocean only sufficiently far from the coastline.
We refer to \cite{lotto2015tsunami} and references therein for details on the model derivation.

The model is formulated as a mixed problem in terms of the vector-valued velocity field $u$, the scalar-valued pressure field $p$, and the surface gravity wave height $\eta$:
\be
   \left\{
   \begin{aligned}
      \rho \frac{\p u}{\p t} + \nabla p &= 0, & \Omega \times (0,T), \\
      K^{-1} \frac{\p p}{\p t} + \nabla \cdot u &= 0, & \Omega \times (0,T), \\
      p &= \rho g \eta, & \p \Omega_{\text s} \times (0,T), \\
      \frac{\p \eta}{\p t} &= u \cdot n, & \p \Omega_{\text s} \times (0,T), \\
      u \cdot n &= - \frac{\p b}{\p t}, & \p \Omega_{\text b} \times (0,T), \\
      u \cdot n &= Z^{-1} p, & \p \Omega_{\text a} \times (0,T), \\
   \end{aligned}
   \right.\label{eq:fwd-pde}
\ee
with homogeneous initial conditions; $K$ and $\rho$ are respectively the bulk modulus and density of seawater, $Z = \rho c$ is the acoustic wave impedance, $c = \sqrt{K/\rho}$ is the speed of sound in seawater, and $\p b / \p t$ is the seafloor motion; the spatial domain is denoted by $\Omega$ with boundaries $\p \Omega_{\text s}$ (sea surface), $\p \Omega_{\text b}$ (sea bottom), and $\p \Omega_{\text a}$ (lateral, absorbing boundaries); $n$ is the outward unit normal; and the temporal domain is $(0,T)$.

The acoustic--gravity model \eq{fwd-pde} represents an LTI dynamical system of the form \eq{LTI} and fits into the Bayesian inference framework laid out in the preceding section.
The model parameter $m$ is given by the spatiotemporal seafloor motion, here defined as the seafloor inward normal velocity, $m(x,t) \coloneqq \p b(x,t) / \p t$, $(x,t) \in \p \Omega_{\text b} \times (0,T)$.
The state is comprised of the velocity, pressure, and surface gravity wave height unknowns, $w \coloneqq (u, p, \eta)$, with homogeneous initial state.
The observables of the system are pointwise observations of the pressure field, assumed to be obtained from pressure sensors mounted on the sea bottom, i.e.~$d = \mc B w = p(x^d_j, t^d_i)$ where $x^d_j \in \p \Omega_{\text b}$, $j = 1,\ldots,\numdata$, are the seafloor sensor locations and $t^d_i \in (0,T)$, $i = 1,\ldots,\numtime^d$, are the time instances at which observations are made.
The model QoIs are observations of the surface wave height, i.e.~$q = \mc B_q w = \eta(x^q_j, t^q_i)$ where $x^q_j \in \p \Omega_{\text s}$, $j = 1,\ldots,\numqoi$, and $t^q_i \in (0,T)$, $i = 1,\ldots,\numtime^q$, are specific locations and time instances pertinent to tsunami early warning.

The forward problem \eq{fwd-pde} can be cast into a mixed variational formulation and then discretized with finite elements and explicit time-stepping (see \cite{henneking2025tsunami} for details).
The discrete expressions of the p2o map, $\datavec = \ptomat \paramvec$ (cf.~\eq{p2o}), and the p2q map, $\qoivec = \ptqmat \paramvec$ (cf.~\eq{p2q}), are then clearly defined as the mappings of the parameter input $\paramvec$ (seafloor motion) to the observables $\datavec$ (sea bottom pressure) and the QoIs $\qoivec$ (surface gravity wave height) via the solution to the discretized equations of the acoustic--gravity PDE model \eq{fwd-pde}.

The inverse problem can then be stated as follows:
given pressure recordings $\obsvec$ from sensors on the seafloor, infer the spatiotemporal seafloor motion $\paramvec$ in the subduction zone.
And the tsunami prediction problem is:
given inferred parameters $\paramvec$, forward predict the surface gravity wave heights $\qoivec$.
In the goal-oriented setting, the task is to directly predict $\qoivec$ from the observations $\obsvec$ without the need to explicitly infer $\paramvec$ (even if this inverse problem is embedded within the goal-oriented inference of $\qoivec$).

\paragraph{Adjoint equations}
Solving the goal-oriented inverse problem (i.e.~characterizing the Bayesian posterior \eq{qoi-posterior}) requires us to have access to the Hessian action \eq{Hessian} which involves the adjoints of the p2o and p2q maps, $\ptomat^*$ and $\ptqmat^*$.
There are two common approaches for deriving the adjoint problems: discretize-then-optimize (DTO) and optimize-then-discretize (OTD).
Each approach has its benefits and drawbacks, see \cite{gunzburger2002control} for a detailed discussion.
While the discretized adjoint problems derived by DTO and OTD may
agree under certain circumstances, in general 
they need not agree, or else they may agree only in the discretization
limit.\footnote{Depending on the discretization scheme, DTO and OTD may yield consistent adjoint problems. For the numerical example presented in Section~\ref{sec:numerical}, this is in fact the case. An adjoint-consistent discretization ensures correctness of gradients.}
If we had direct access to the discretized maps $\ptomat$ and $\ptqmat$, then the DTO approach could be realized by simply transposing these maps. Typically, this is not a practical method because we usually are not able to precompute and store these matrices but only have access to their action on a vector which requires solving the PDEs.\footnote{In the case where $\ptomat$ and $\ptqmat$ represent discretized LTI systems, they can be compactly stored; however, precomputing them efficiently still requires access to the action of the adjoints $\ptomat^*$ and $\ptqmat^*$ (as will be shown in Section~\ref{sec:methodology}).} 
Here, we present the adjoint equations at the infinite-dimensional level, consistent with an OTD approach.
Similar to the forward problem \eq{fwd-pde}, the adjoint PDEs are given as a mixed problem in terms of the adjoint velocity field $\tau$, the adjoint pressure $v$, and the adjoint surface gravity wave height $\xi$:

\be
   \left\{
   \begin{aligned}
      \rho \frac{\p \tau}{\p t} + \nabla v &= 0, & \Omega \times (0,T), \\
      K^{-1} \frac{\p v}{\p t} + \nabla \cdot \tau &= 0, & \Omega \times (0,T), \\
   \end{aligned}
   \right.\label{eq:adj-pde}
\ee
with homogeneous terminal conditions; the BCs for the adjoint p2o map are
\be
   \left\{
   \begin{aligned}
      v &= \rho g \xi, & \p \Omega_{\text s} \times (0,T), \\
      \frac{\p \xi}{\p t} &= \tau \cdot n, & \p \Omega_{\text s} \times (0,T), \\
      \tau \cdot n &= - \mc B^* \tilde d, & \p \Omega_{\text b} \times (0,T), \\
      \tau \cdot n &= -Z^{-1} v, & \p \Omega_{\text a} \times (0,T), \\
   \end{aligned}
   \right.\label{eq:adj-p2o-bc}
\ee
where $\tilde d$ is the source of the adjoint problem, defined at the observation points.\footnote{Technically, the adjoint observation operator $\mc B^*$ maps to the adjoint state space, so $\mc B^* \tilde d$ could be written as a source term on the right-hand side of \eq{adj-pde}.}

For the adjoint p2q map, the BCs change; in particular, the boundary source term, denoted by $\tilde q$, is now defined on the sea surface, and we instead obtain homogeneous BCs on the sea bottom:
\be
   \left\{
   \begin{aligned}
      v &= \rho g \xi, & \p \Omega_{\text s} \times (0,T), \\
      \frac{\p \xi}{\p t} &= \tau \cdot n + \mc B_q^* \tilde q, & \p \Omega_{\text s} \times (0,T), \\
      \tau \cdot n &= 0, & \p \Omega_{\text b} \times (0,T), \\
      \tau \cdot n &= -Z^{-1} v, & \p \Omega_{\text a} \times (0,T). \\
   \end{aligned}
   \right.\label{eq:adj-p2q-bc}
\ee

Analogous to the forward problem, these adjoint problems can be cast into mixed variational formulations and then discretized with finite elements and a time-stepping scheme to obtain discrete representations of the adjoint p2o and p2q maps $\ptomat^*$ and $\ptqmat^*$.
Each action of $\ptomat^*$ and $\ptqmat^*$ then corresponds to mapping the respective source terms $\tilde \datavec$ and $\tilde \qoivec$ to the parameter space via the solution to the discretized equations of the adjoint acoustic--gravity PDEs \eq{adj-pde}+\eq{adj-p2o-bc} and \eq{adj-pde}+\eq{adj-p2q-bc}, respectively.

\paragraph{Hessian action}
The linear inverse problem \eq{mmap-problem} can be solved with an iterative solver like the conjugate gradient method,
which requires numerous actions of the Hessian on a vector.
Explicit construction of the discrete Hessian itself is usually prohibitive, because it necessitates as many forward PDE solves as there are parameters, leading to the need for performing matrix-free Hessian matvecs instead.
Each action of the adjoint-based Hessian $\vb H$ involves matvecs with both $\ptomat$ and $\ptomat^*$ (cf.\ \eq{Hessian}).
Matrix-free Hessian matvecs therefore require solution of a pair of forward/adjoint PDE solves per Hessian action.
Unfortunately, solving \eq{mmap-problem} in this way is not feasible in the case of the tsunami inverse problem.
Its massive scale and the need for real-time solution preclude solving the forward/adjoint PDE model at the time of inference.
Instead, as will be shown in Section~\ref{sec:methodology}, we exploit the block Toeplitz structure of $\ptomat$ and $\ptomat^*$, as well as $\ptqmat$ and $\ptqmat^*$, which implies that we can precompute and store these matrices, and it enables Hessian matvecs that are orders of magnitude faster and more efficient than matrix-free actions via pairs of forward/adjoint PDE solves.

\paragraph{Discretization of the CSZ tsunami inverse problem}
In order to present the reader with concrete examples of the memory and computational complexities involved in our real-time inference and prediction methodology (presented in Section~\ref{sec:methodology}), we provide estimates based on the problem size for the full-scale CSZ digital twin for tsunami early warning.
We re-emphasize that the methodology presented below can be applied to many other applications whose problem dimensions may be quite different from these estimates.

The CSZ is more than 1000 km long, stretching from off the coast of Northern California to beyond the Northwest peak of Vancouver Island.
Discretizing the full-order acoustic--gravity PDE model along the entire subduction zone and extending near the Pacific--Northwest coastline yields a forward problem with $\numstate \sim 10^{10}$ states and $\numparam \sim 10^{6}$ parameters per time step (i.e.~$10^{10}$ and $10^6$ spatial degrees of freedom, respectively).\footnote{The number of time steps for the state is determined by the Courant--Friedrichs--Lewy (CFL) condition.}
Sea bottom sensors in the Cascadia region are very sparse but various efforts are underway to instrument the CSZ, as well as other subduction zones around the world, to provide real-time data for tsunami early warning.
For the CSZ digital twin, observations are assumed to be obtained from seafloor pressure measurements recorded by hundreds of sensors, i.e.~$\numdata \sim 10^2$, and tsunami predictions are made in a few specified forecasting locations, i.e.~$\numqoi \sim 10$ (e.g.~near major population centers).
The dominant frequencies of earthquake-induced seafloor motion and the ensuing ocean acoustic waves are typically below 1~Hz (though sea bottom pressure sensors may record higher frequencies)~\cite{li2009ocean}.
Given the proximity of the CSZ to the coastline, where the tsunami may reach population centers within 20 minutes of the earthquake, the observation time prior to inference and prediction with the CSZ digital twin to provide early warning is on the order of a few hundred or at most one thousand seconds, so $\numtime \sim 10^3$.

\section{Methodology}
\label{sec:methodology}

The Bayesian posteriors of the parameters $\paramvec$ and the QoIs $\qoivec$, derived in Section~\ref{sec:background}, are respectively given by $\pipost(\paramvec | \obsvec) \propto \mc N(\bfmmap, \bfpostcovm)$ and $\pipost(\qoivec | \obsvec) \propto \mc N \left( \bfqmap, \bfpostcovq \right)$, where
\begin{align}
    &\bfpostcovm = \left(\ptomat^* \bfinvnoise \ptomat + \bfinvpriorcovm \right)^{-1}
    &&,\quad \bfmmap = \bfpostcov \left(\ptomat^* \bfinvnoise \obsvec + \bfinvpriorcovm \bfmprior \right) ,
    \label{eq:map-point} \\
    &\bfpostcovq = \ptqmat \bfpostcovm \ptqmat^*
    &&,\quad \bfqmap = \ptqmat \bfmmap .
    \label{eq:qoi-map-point}
\end{align}

The ultimate goal of Bayesian inference is to characterize these posteriors; for example, computing their means, estimating their covariances, or drawing samples. For inverse problems where the effective rank of the Hessian $\vb{H}\coloneqq\bfinvpostcovm$ is high (e.g.~hyperbolic systems like the acoustic--gravity wave equations), it takes many iterations for an iterative solver such as the conjugate gradient method to solve for $\vb{m}_{\text{map}}$ in~\eq{map-point} using actions of $\vb H$. Moreover, conventionally each application of $\vb{H}$ requires a pair of forward and adjoint wave equation solves, which is expensive and precludes real-time inference. Here, we present an alternative: the tsunami dynamics are governed by an autonomous system whose operator does not have explicit time dependence. This fact (and causality) manifests in the time-shift invariance of the p2o and p2q maps ($\ptomat$ and $\ptqmat$), which results in $\ptomat$ and $\ptqmat$ being block lower-triangular Toeplitz matrices (cf.~\eq{ShiftInvariance}), enabling extremely fast and efficient matvecs with $\ptomat$, $\ptqmat$, and their respective adjoints (see below). In this section, we introduce a 4-phase framework for real-time inference and prediction; this framework relies on the following key ideas:
\begin{itemize}
   \item efficient and fast FFT-based multi-GPU-accelerated actions of the discretized p2o and p2q operators and their respective adjoints;
   \item compact storage of the p2o and p2q maps in block Toeplitz form;
   \item an offline--online decomposition of the inference and prediction problem into precomputation phases (offline) and a real-time inference and prediction phase (online); and
   \item a shift of the computation of the inverse operator from the high-dimensional parameter space to the much lower-dimensional data space.
\end{itemize}

In the remainder of this section, $\numtime$ refers to the number of time steps in the discrete parameter field, data, and QoIs;\footnote{For ease of presentation, without loss of generality, we use the same time step for parameters, observables, and QoIs. In practice, each one may have a different time step as needed (see numerical example in Section~\ref{sec:numerical}).} $\numparam$ refers to the number of spatial degrees of freedom in the discrete parameter field; $\numdata$ to the number of sensor locations; and $\numqoi$ to the number of spatial locations for QoI prediction.
While our methodology does not assume specific values for any of these parameters, an efficient realization of the framework does require the following assumptions:
\be
   \numdata \ll \numparam,\
   \numqoi \ll \numparam,\
   \text{and } \numtime \gg 1.
\ee

In the context of inverse problems governed by autonomous dynamical systems, these assumptions are not very restrictive.
In particular, many applications are data-constrained, i.e.~they seek to infer high-dimensional parameter fields from sparse measurements.
And in the goal-oriented setting, we often aim to predict a limited number of QoIs which are typically much fewer than the parameter degrees of freedom.
Recall from the preceding section that our estimates for the full-scale CSZ digital twin are $\numparam \sim \! 10^6$, $\numdata \sim \! 10^2$, $\numqoi \sim \! 10$, $\numtime \sim \! 10^3$.
It is easily verified that these problem dimensions fit well within the intended scope of the framework.

\paragraph{FFT-based GPU-accelerated p2o and p2q matvecs}

To solve the parameter inference and QoI prediction problem (i.e.~characterize the corresponding Bayesian posteriors), we require many matvecs with the p2o and p2q maps ($\ptomat$ and $\ptqmat$) and their respective adjoints ($\ptomat^*$ and $\ptqmat^*$). Matrix-free actions of those maps and their adjoints on a vector involve the solution of the forward and adjoint acoustic--gravity wave equations, respectively, making them computationally expensive.

However, if we have direct access to the matrices $\ptomat$ and $\ptqmat$, we can exploit their block lower-triangular Toeplitz structure \eq{ShiftInvariance}.
Fast matvecs are achieved by embedding the respective block Toeplitz matrix within a block circulant matrix, which is block-diagonalized by the discrete Fourier transform (DFT) \cite{gray2006toeplitz}.
The DFT only needs to be done once for each matrix.
The matvec then becomes a block-diagonal matvec operation in Fourier space.\footnote{The matvec in Fourier space can also be computed as an elementwise vector operation for each block followed by a row-wise reduction over blocks. However, this requires strided memory access patterns which negatively impacts GPU performance; see \cite{venkat2024fft} for details.}
Moreover, because the FFT is a unitary operator, the action of the adjoint p2o and p2q maps correspond to simply applying the conjugate-transposed diagonal blocks in Fourier space, eliminating the need to separately store the Fourier-transformed forward and adjoint maps.
These FFT-based p2o and p2q matvecs can be implemented efficiently on multi-GPU clusters~\cite{venkat2024fft}.

In comparison to the conventional matrix-free actions involving forward and adjoint PDE solves, the FFT-based algorithm trades off computational complexity (PDE solves) with memory complexity (storing the matrices).
Conventionally, computing or storing the full p2o map is not tractable for large-scale inverse problems, because doing so requires one forward PDE solve per column of the matrix, and the number of columns corresponds to the dimension of the parameter space ($\sim \! 10^9$). For linear inverse problems governed by autonomous dynamical systems, however, the autonomous structure of the p2o and p2q maps both reduces the number of PDE solves required and enables compact storage of these maps (see Phase~1 below).

The main advantages of the FFT-based matvec algorithm are that it (1) avoids PDE solves entirely, thus making it independent of the cost of state discretization, time-stepping constraints from CFL condition, and computing on unstructured mesh-based data structures dictated by the PDE discretization; and (2) relies on arithmetic operations that involve data that is contiguous in memory and leverages routines from libraries such as cuBLAS and rocBLAS that are well known to achieve high performance on GPUs. 
A quantitative comparison between the FFT-based GPU-accelerated algorithm and the conventional approach of matrix-free actions via forward and adjoint PDE solves for a representative test problem is presented in Section~\ref{sec:numerical}.

\paragraph{Phase 1 (offline): Precompute p2o and p2q maps}

Phase 1 of our framework is concerned with precomputing and storing the p2o and p2q maps, so that their actions are later available (during Phases 2--4) at the much lower computational cost enabled by the FFT-based algorithm and with the speedup provided due to the algorithm's suitability for multi-GPU execution.

The p2o and p2q maps can be computed either from forward PDE solutions (one per column) or from adjoint PDE solutions (one per row). Because the data dimension ($\numdata \numtime \sim \! 10^5$) and QoI dimension ($\numqoi \numtime \sim \! 10^4$) are each much smaller than the parameter dimension ($\numdata \numtime \sim \! 10^9$), computing rows via adjoint PDE solutions is the more efficient approach.
The time-shift invariance of the discretized p2o and p2q maps implies that only one adjoint PDE solve per sensor location is required to construct $\ptomat$, and one adjoint PDE solve per QoI spatial location to construct $\ptqmat$.
For $\ptomat$, these adjoint PDEs are solved backward in time with pressure ``sources'' at each sensor location (cf.\ $\tilde d$ in \eq{adj-p2o-bc}); for $\ptqmat$, they are solved with pressure sources at each QoI point in space (cf.\ $\tilde q$ in \eq{adj-p2q-bc}). Translating these solutions in time gives the entire $\ptomat$ and $\ptqmat$. 
The block Toeplitz structure of $\ptomat$ and $\ptqmat$ also reduces the memory requirements to store these matrices by a factor of $\numtime$.
For the CSZ digital twin, this means only $\sim \! 0.8$~TB memory is required to store $\ptomat$ instead of $\sim \! 800$~TB, and $\sim \! 80$~GB to store $\ptqmat$ instead of $\sim \! 80$~TB.

While $\numdata + \numqoi \sim \! \mc{O}(10^2)$ adjoint acoustic--gravity wave-equation solves incur considerable computational expense, this entire computation happens offline and only needs to be done once.
For the discretized CSZ digital twin (with $\sim \! 10^{10}$ states and $\sim \! 10^3$ time steps), a scalable explicit finite element solver executed on a moderately large compute cluster is able to compute one adjoint acoustic--gravity wave-equation solution per hour.\footnote{This time estimate was obtained from a benchmark of the acoustic--gravity wave equation solve performed on 512 compute nodes of TACC's \emph{Frontera} supercomputer using an explicit finite element solver implemented in MFEM.}

\paragraph{Phase 2 (offline): Compute a compact representation of $\bfpostcovm$}

For the CSZ digital twin, the discretized parameter dimension is $\numparam \numtime \sim \! 10^{9}$.
The linearity of the inverse problem enables us, in principle, to precompute both the mapping from observed data $\obsvec$ to the MAP point $\bfmmap$ as well as the posterior covariance $\bfpostcovm$:
\begin{align}
   \bfmmap &=
   \underbrace{\vb H^{-1} \ptomat^* \bfinvnoise}_{\text{precompute}} \obsvec
   +
   \underbrace{ \vb H^{-1} \bfinvpriorcovm \bfmprior}_{\text{precompute}} \, , 
   \label{eq:map-point-precompute} \\
   \bfpostcovm &= \vb H^{-1} =
   \underbrace{
   \left( \ptomat^* \bfinvnoise \ptomat + \bfinvpriorcovm \right)^{-1}}_{\text{precompute}} \, .
\end{align}
Unlike the p2o map, the inverse of the Hessian does not inherit the block Toeplitz structure from the discretized autonomous dynamical system, so computing and storing the full covariance matrix $\bfpostcovm$ of size $\sim \! 10^9 \times 10^9$ is intractable.
Moreover, the forward model~\eq{fwd-pde} is hyperbolic, implying a relatively slow decay of the eigenvalues of the Hessian (see Section~\ref{sec:numerical}), thereby prohibiting the effectiveness of traditional low-rank decomposition methods~\cite{ghattas2021learning}.

To circumvent this issue, we employ the Sherman--Morrison--Woodbury formula~\cite{sherman1950, woodbury1950}, with which $\bfpostcovm$ can be rewritten as
\be
   \bfpostcovm = \bfpriorcovm \left( \Imat - \ptomat^* \Kmat^{-1} \ptomat \bfpriorcovm \right) ,
   \label{eq:woodbury}
\ee
where $\Kmat = \bfnoise + \ptomat \bfpriorcovm \ptomat^*$.
Expressing $\bfpostcovm$ in this form shifts the problem of inverting an operator of the dimension of the parameter space $\numparam \numtime$ to one of the dimension of the data space $\numdata \numtime$. For the CSZ digital twin, this corresponds to a matrix $\Kmat$ of size $10^5 \times 10^5$, whereas the original $\bfpostcovm$ was of size $10^{9} \times 10^{9}$. 
This matrix can be Cholesky-factorized offline, and triangular solves are done in the online phase.

To compute matrix $\Kmat$, Phase~2 of our framework relies on the efficient FFT-based algorithm to perform actions of the p2o map $\ptomat$ and its adjoint $\ptomat^*$ precomputed in Phase~1.
Computing $\Kmat$ column-by-column requires $10^5$ matvecs with $\ptomat$ and $\ptomat^*$ each.
Using the conventional approach of matrix-free actions via forward and adjoint PDE solves to perform $10^5$ matvecs would have made this task infeasible.
As we will see for a numerical example in Section~\ref{sec:numerical}, the FFT-based algorithm can perform these matvecs for large-scale problems within a reasonably short time with moderate compute resources.

Computing $\Kmat$ also requires actions of the noise covariance $\bfnoise$ and the prior covariance $\bfpriorcovm$. 
We assume that $\bfnoise$ is a diagonal matrix, making the cost of applying $\bfnoise$ or $\bfnoise^{-1}$ to a vector negligible.
For the prior covariance, a common choice is to employ a smoothing prior given in the form of a discretized inverse elliptic differential operator \cite{ghattas2021learning, stuart2010inverse, lindgren2011explicit}.
Performing $10^5$ solves of a discretized elliptic PDE operator of size $10^9$ is a tractable problem thanks to the existence of fast and scalable solvers for such systems, but it can become a bottleneck in the computation of $\Kmat$.
However, many times the elliptic PDE operator whose inverse generates the prior covariance may itself be time-invariant (we use one such example in Section~\ref{sec:numerical}), so that the matrix $\Gmat^* \coloneqq \bfpriorcovm \ptomat^*$ is also block Toeplitz.
In this case, $\Gmat^*$ can be precomputed at the cost of only $\numdata \sim \! 10^2$ prior solves (which can be batched over time steps) and stored at the same expense as storing $\ptomat$. Matvecs with $\Gmat^*$ can then be performed with the FFT-based algorithm and incur the same computational cost as matvecs with $\ptomat^*$.
With the definition of $\Gmat^*$, the updated expression for the Hessian inverse is
\be
   \bfpostcovm =
   \left( \Imat - \Gmat^* \Kmat^{-1} \ptomat \right) \bfpriorcovm ,
   \label{eq:fast-Hessian-matvec}
\ee
where $\Kmat = \bfnoise + \ptomat \Gmat^*$.

It must be emphasized that Hessian actions given by \eq{fast-Hessian-matvec} are exact up to rounding error.
In other words, compact representation of $\bfpostcovm$ here only refers to an efficient representation of the inverse operator in the data space dimension, which was obtained through applying the Sherman--Morrison--Woodbury formula to $\bfpostcovm$ (cf.~\eq{woodbury}), and it does \emph{not} refer to an approximation of the Hessian action through representing $\bfpostcovm$ in some lower-dimensional space.

Having precomputed $\Gmat^*$ and prefactorized $\Kmat$, Hessian actions via \eq{fast-Hessian-matvec} with the multi-GPU-accelerated FFT-based matvecs of $\ptomat$ and $\Gmat^*$ are extremely fast and efficient, and they enable the real-time inference of the parameter MAP point, as will be shown in Phase~4.
In fact, as we will demonstrate with numerical results in Section~\ref{sec:numerical}, the Hessian matvec computation time is now primarily limited by the prior solve of $\bfpriorcovm$.
Despite this efficiency of Hessian matvecs, computation of the full posterior covariance may still be unattainable due to the memory cost of storing columns of $\bfpostcovm$ as well as the computational cost of performing as many as $\sim \! 10^9$ Hessian matvecs.
However, obtaining a subset of the parameter uncertainties from the pointwise variance field, given by $\diag(\bfpostcovm)$, is now feasible by computing columns of $\bfpostcovm$.

Lastly, if the prior mean is non-zero ($\bfmprior \ne \vb{0}$), we can now easily precompute its contribution to the MAP point of the posterior (cf.~\eq{map-point-precompute}):
\be
   \bfmmapprior \coloneqq
   \bfpostcovm \bfinvpriorcovm \bfmprior =
   \left( \Imat - \Gmat^* \Kmat^{-1} \ptomat \right) \bfmprior .
   \label{eq:map-point-prior}
\ee

\paragraph{Phase 3 (offline): Compute QoI uncertainties and data-to-QoI map}

Many applications of interest, including tsunami early warning, ultimately require a real-time forecast (with uncertainties) for one or multiple specified QoIs.
To achieve this, we can exploit the fast Hessian actions enabled by the compact representation of $\bfpostcovm$, precomputed in Phase~2, to compute uncertainties for the QoIs $\qoivec$ via
\be
   \bfpostcovq 
   = \ptqmat \bfpostcovm \ptqmat^*
   = \ptqmat \left( \Imat - \Gmat^* \Kmat^{-1} \ptomat \right) \bfpriorcovm \ptqmat^* .
   \label{eq:qoi-covariance}
\ee
In the case of a time-invariant prior, we can use the same trick as in computing $\Gmat^*$ (for faster Hessian matvecs) to accelerate computation of \eq{qoi-covariance} by precomputing the block Toeplitz matrix $\Gmat_q^* \coloneqq \bfpriorcovm \ptqmat^*$ via $\numqoi \sim \! 10$ prior solves (one for each column of the last block column of $\ptqmat^*$).
The posterior covariance for the QoIs can then be written as
\be
   \bfpostcovq 
   = \ptqmat \left( \Imat - \Gmat^* \Kmat^{-1} \ptomat \right) \Gmat_q^* .
   \label{eq:fast-qoi-covariance}
\ee
The covariance matrix $\bfpostcovq$ of size $\sim \! 10^4 \times 10^4$ can now be computed through $\sim \! 10^4$ actions of $\bfpostcovq$ via \eq{fast-qoi-covariance}, giving access to the Bayesian posterior of $\qoivec$ once $\bfqmap$ has been computed.

During the online phase of real-time inference and prediction, described in Phase~4, the QoI MAP point $\bfqmap$ could, in principle, be computed through push-forward of the inferred parameters $\bfmmap$ via the operator $\ptqmat$ (cf.~\eq{qoi-map-point}).
Using the FFT-based matvec of $\ptqmat$, this is entirely feasible in real time.
However, this requires (1) first solving the real-time inference of $\bfmmap$; (2) storing the matrices $\Kmat$, $\ptomat$, $\ptqmat$, and $\Gmat^*$ for use during the online phase; and (3) online access to a multi-GPU cluster to compute fast actions of these operators.

A more efficient alternative for the real-time QoI forecast is available by precomputing the data-to-QoI (d2q) map
\be
   \Qmat \coloneqq \ptqmat \bfpostcovm \ptomat^* \bfinvnoise
   = \ptqmat \left( \Imat - \Gmat^* \Kmat^{-1} \ptomat \right) \Gmat^* \bfinvnoise ,
   \label{eq:d2q-map}
\ee
which can then be used to directly map observed data to the QoI without the need to first infer the parameter $\bfmmap$.
This in turn drastically reduces the storage and compute requirements for the QoI prediction during the online phase (see Phase~4 below).
The d2q map is of size $\sim \! 10^4 \times 10^5$, so it can be computed through $\sim \! 10^4$ actions of $\Qmat$ via \eq{d2q-map}.

We remark that a similar computation could have been performed to precompute the data-to-parameter map given by $\bfpostcovm \ptomat^* \bfinvnoise$. However, for the problem scale we are interested in, the memory requirements for storing this mapping ($\sim \! 800$~TB) together with the I/O time required for reading it from disk during the online phase make this difficult to use in real time. In other applications---particularly where the parameter dimension is a few orders of magnitude smaller than for the CSZ digital twin---this may be an attractive approach.

Finally, the contribution of a non-zero prior mean ($\bfmprior \ne \vb{0}$) to the MAP point of the QoI can also be precomputed (cf.~\eq{qoi-map-point} and \eq{map-point-prior}):
\be
   \bfqmapprior \coloneqq
   \ptqmat \bfmmapprior =
   \ptqmat \bfpostcovm \bfinvpriorcovm \bfmprior = 
   \ptqmat \left( \Imat - \Gmat^* \Kmat^{-1} \ptomat \right) \bfmprior .
\ee

\paragraph{Phase 4 (online): Real-time inference and prediction}

After the various stages of offline computation during Phases~1--3, the online phase reduces to the tasks of (1) real-time inference of the parameters $\bfmmap$ and (2) real-time prediction of the QoIs $\bfqmap$.

Using the derived expression \eq{fast-Hessian-matvec} for fast Hessian matvecs, the inference of $\bfmmap$ from real-time observations $\obsvec$ can be calculated via
\be
   \bfmmap = 
   \left( \Imat - \Gmat^* \Kmat^{-1} \ptomat \right) \Gmat^* \bfinvnoise \obsvec
   + \bfmmapprior
   \label{eq:real-time-parameter}
\ee
in real time (order of seconds) for $\sim \! 10^9$ parameters.
This real-time inference relies on the FFT-based GPU-accelerated matvecs with the precomputed block Toeplitz matrices $\ptomat$ and $\Gmat^*$, and triangular solves of the prefactorized matrix $\Kmat$.

The real-time prediction for the QoI employs the precomputed d2q map $\Qmat$ (a dense matrix of size $\sim \! 10^4 \times 10^5$). The prediction of $\bfqmap$ can then be obtained directly from observations $\obsvec$ via
\be
   \bfqmap = \Qmat \obsvec + \bfqmapprior
   \label{eq:real-time-qoi}
\ee
for $\sim \! 10^4$ QoIs within fractions of a second.
Importantly, due to its relatively small size, the d2q map can be stored with about $\sim \! 8$~GB of memory, and the dense matvec of $\Qmat$ is computationally cheap. 
This makes the action of $\Qmat$ available for online computation on a laptop.
With precomputed $\Qmat$, this QoI prediction is accomplished without the need for inference of $\bfmmap$; therefore, \eq{real-time-qoi} enables high-fidelity (exact up to rounding error) goal-oriented Bayesian inference in real time without the need for online access to large-scale compute resources.

If this goal-oriented setting allows us to obtain real-time predictions of the QoIs directly from the observed data, then one might ask why we still care for the inference of the parameters.
Inferring the earthquake-induced spatiotemporal seafloor motion is of interest, because this inverse solution can be used as a source for tsunami propagation in other models, including nonlinear tsunami propagation models that account for coastal inundation and other features that are not represented in the linearized acoustic--gravity PDE model.

\section{Numerical Example}
\label{sec:numerical}

In this section, the end-to-end digital twin framework---from noisy observations via parameter inference to QoI prediction with uncertainties---is demonstrated with a 3D test configuration of the acoustic--gravity wave equations. We employ synthetic observations (obtained from the PDE model) with additive noise as pressure recordings from seafloor sensors to infer the underlying spatiotemporal seafloor motion and forward predict tsunami wave heights at specified locations in real time.

For the Gaussian prior, we use prior mean $\bfmprior = \vb{0}$, and the prior covariance $\bfpriorcovm$ is generated from the inverse of a discretized elliptic differential operator:
\be
   \bfpriorcovm \coloneqq \left( 
   \alpha_1 I_h - 
   \alpha_2 \Delta_h
   \right)^{-2} ,
   \label{eq:prior-cov}
\ee
where $\alpha_1 > 0$, $\alpha_2 > 0$, $I$ is identity, $\Delta$ denotes the Laplacian $(\p^2 / \p x^2 + \p^2 / \p y^2)$, and $h$ indicates discretization;
we employ Robin boundary conditions to mitigate boundary artifacts~\cite{daon2018mitigating}.
This choice of prior is equivalent to a Mat\'ern covariance kernel~\cite{lindgren2011explicit} and it guarantees bounded pointwise variance~\cite{stuart2010inverse}.
The prior covariance \eq{prior-cov} smoothes rough components of the parameter field.
The premise for using a smoothing prior operator is that the inference of rough components is unstable because they are not sufficiently informed by the data~\cite{ghattas2021learning}.
Note that this prior has time-invariant parameters $\alpha_1$ and $\alpha_2$, and it uses spatial correlation only, leading to a shift-invariant prior covariance matrix with block-diagonal structure.\footnote{The shift invariance of the prior covariance $\bfpriorcovm$ enables efficient matvecs with $\Gmat^* \coloneqq \bfpriorcovm \ptomat^*$ and $\Gmat_q^* \coloneqq \bfpriorcovm \ptqmat^*$ by precomputing and storing $\Gmat^*$ and $\Gmat_q^*$ as block Toeplitz matrices. Without this structure of $\bfpriorcovm$, the number of prior solves is significantly larger during Phases 2--3, making these offline computations more costly.}

\paragraph{Test configuration}

\begin{figure}[htb]
   \includegraphics[width=\textwidth,trim={0 0 0 0pt}]
   {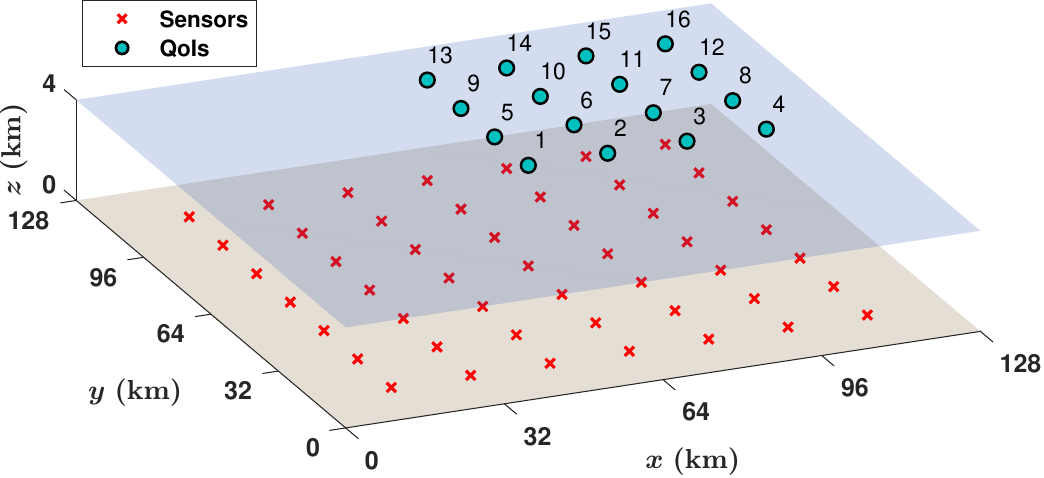}
   \caption{3D spatial domain of a test configuration with the locations of 49 seafloor pressure sensors and 16 predicted QoIs (tsunami wave heights).}
   \label{fig:domain}
\end{figure}

Figure~\fig{domain} illustrates the spatial domain and the locations of sensors and QoIs for the test configuration.
The spatial domain is given by $\Omega \coloneqq (0,128) \text{ km} \times (0,128) \text{ km} \times (0,4) \text{ km}$, and it is discretized on a spatial grid with resolution of $0.25 \text{ km} \times 0.25 \text{ km} \times 0.125 \text{ km}$.
Each grid point on the seafloor has one parameter degree of freedom, leading to a spatial parameter dimension of $\numparam = 513 \times 513 = 263\,169$.
The total number of sensors is $\numdata = 49$; the sensors are placed uniformly on a $7 \times 7$ grid:
\be
   \begin{split}
   \quad x_d &\in \{ 0.125, 0.25, 0.375, 0.5, 0.625, 0.75, 0.875 \} 
   (x_{\max} - x_{\min}) , \\
   \quad y_d &\in \{ 0.125, 0.25, 0.375, 0.5, 0.625, 0.75, 0.875 \} 
   (y_{\max} - y_{\min}) ,
   \end{split}
\ee
where $x_{\min} = y_{\min} = z_{\min} = 0$, $x_{\max} = y_{\max} = 128$, $z_{\max} = 4$, and $(x_d, y_d, z_{\min})$ are the sensor coordinates on the seafloor.
The QoI is the surface-gravity wave height at $\numqoi = 16$ spatial locations with coordinates $(x_q, y_q, z_{\max})$ on the sea surface given by
\be
   \begin{split}
   \quad x_q &\in \{ 0.5, 0.625, 0.75, 0.875 \} 
   (x_{\max} - x_{\min}) , \\
   \quad y_q &\in \{ 0.5, 0.625, 0.75, 0.875 \} 
   (y_{\max} - y_{\min}) .
   \end{split}
\ee
The time domain is given by $(0, T) \coloneqq (0, 50)$~s.
Parameters and data are discretized at the same temporal frequency of 10~Hz, leading to $\numtime^m = \numtime^d =: \numtime = 500$.
The QoI is discretized at a lower temporal resolution of 1~Hz, so that $\numtime^q = 50$.
The total (spatiotemporal) dimensions are respectively $\numparam \numtime = 131\,584\,500$ parameters, $\numdata \numtime = 24\,500$ data, and $\numqoi \numtime^q = 800$ QoIs.

\paragraph{Phase 1 (offline)}

In Phase~1, the adjoint PDEs of the acoustic--gravity wave equations \eq{adj-pde} must be solved repeatedly to precompute the p2o map ($\numdata = 49$ adjoint PDE solves) and the p2q map ($\numqoi = 16$ adjoint PDE solves).
The discrete form of the equations is obtained by discretizing a weak variational formulation of the mixed acoustic--gravity wave problem with Galerkin finite elements in space and explicit 4th-order Runge--Kutta (RK4) time-stepping.
The conforming finite element discretization employs third-order continuous ($H^1$-conforming) elements for the scalar-valued pressure unknown and second-order discontinuous ($L^2$-conforming) elements for the vector-valued velocity unknown (see \cite{henneking2025tsunami} for details).
For the test configuration, this discretization leads to a total state dimension of
$\numstate = 908\,627\,041$ and RK4 time-stepping with $5\,000$ time steps.
The explicit finite element solver for the acoustic--gravity wave equations is implemented in MFEM \cite{anderson2021mfem}.
The implementation is optimized using fast integration via partial assembly with tensor-product hexahedral elements, and it is scalable supporting both CPU and GPU-based computation.

For this numerical example, the adjoint PDEs were solved on the \emph{Stampede3} supercomputer at the Texas Advanced Computing Center (TACC) using 32~SPR compute nodes with a total of $3\,584$~Sapphire Rapids cores (112 cores per node).\footnote{See \url{https://docs.tacc.utexas.edu/hpc/stampede3} for TACC's \emph{Stampede3} system configuration.}
Each one of the 49 + 16 adjoint PDE solves computes and stores one column of the adjoint p2o and p2q map, respectively, each of the size of the parameter dimension ($\sim \! 1.05$~GB), requiring total storage of 52~GB for the p2o map and 17~GB for the p2q map.
The time per solve is approximately eight minutes, leading to a total Phase~1 compute time of $\sim \! 6.5$~hours for the p2o map and $\sim \! 2.1$~hours for the p2q map; however, all of these solves can be done simultaneously to accelerate the computation if needed.

\begin{figure}[htb]
   \centering
   \includegraphics[width=\textwidth,clip,trim={0 0 0 0}]
   {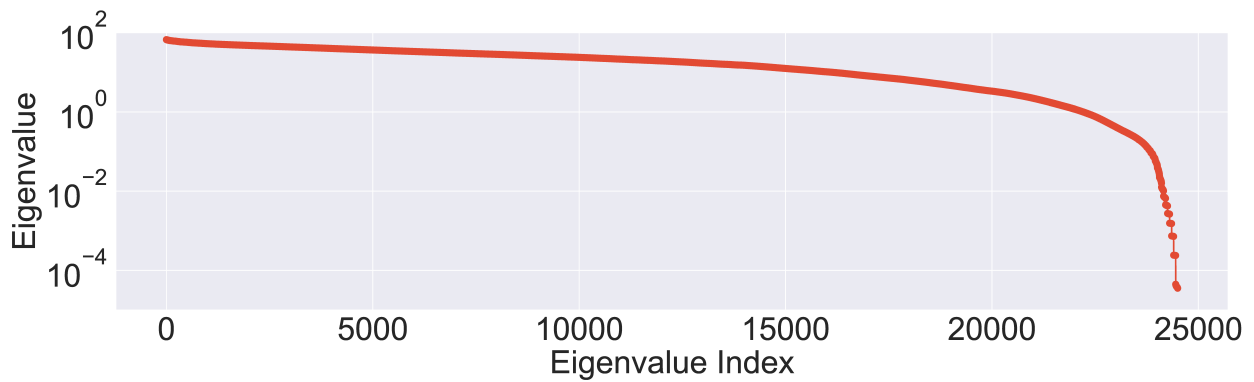}
    \caption{Eigenvalues of $\ptomat \ptomat^*$ for the test configuration. The slow eigenvalue decay is indicative of hyperbolic PDE operators which do not readily admit compact representations.}
   \label{fig:eigenvalues}
\end{figure}

Figure~\fig{eigenvalues} depicts the eigenvalues of $\ptomat \ptomat^*$ for the test configuration.
The plot shows that the eigenvalues decay by only two orders of magnitude within the first 90 percent of the eigenvalues.
This slow eigenvalue decay, indicative of hyperbolic PDE operators, exemplifies the difficulty of finding accurate compact representations of the operator in a lower-dimensional space via low-rank approximation or similar techniques (see our discussion on surrogate modeling in Section~\ref{sec:intro}).

\paragraph{Phase 2 (offline)}

After completion of Phase~1, the finite element solver is no longer required.
Instead of using forward and adjoint PDE solves, the actions of the p2o and p2q maps (and their respective adjoints) are now performed using the multi-GPU-accelerated FFT-based algorithm.
On TACC's \emph{Lonestar6} GPU cluster\footnote{See \url{https://docs.tacc.utexas.edu/hpc/lonestar6} for TACC's \emph{Lonestar6} system configuration.}, with 24 Nvidia A100 GPUs (distributed across eight compute nodes), the runtime reduces to $\sim \! 11$~ms per $\ptomat$ or $\ptomat^*$ matvec and $\sim \! 3$~ms per $\ptqmat$ or $\ptqmat^*$ matvec.
Compared to performing (adjoint) PDE solves with the explicit finite element solver, this corresponds respectively to a $\sim \! 44\,000 \times$ and $\sim \! 160\,000 \times$ speedup for the (adjoint) p2o and (adjoint) p2q matvecs.
These incredibly fast actions of the p2o and p2q maps (and their respective adjoints) empower the remaining computations done offline (Phases~2--3) and the real-time inference and prediction performed online (Phase~4).
All of the reported runtimes in Phases~2--4 are based on the same system configuration with 24~Nvidia A100~GPUs; in some cases, which will be pointed out, computations are done on a single A100~GPU.

Next up is precomputing $\Kmat = \bfnoise + \ptomat \Gmat^*$, where  $\Gmat^* = \bfpriorcovm \ptomat^*$. 
With our choice of prior, $\Gmat^*$ is block Toeplitz, so we only need to perform $\numdata$ prior solves to precompute and store $\Gmat^*$. 
Storage costs of $\Gmat^*$ are equivalent to storing $\ptomat$ ($\sim \! 52$~GB).
The time per prior solve, performed on a single A100~GPU with Nvidia cuDSS\footnote{\url{https://developer.nvidia.com/cudss}}, is $\sim \! 84$~ms. 
Thus, the total time needed to compute $\Gmat^*$ is $\numdata \cdot 84$~ms $\sim \! 4.2$~s.
The matrix $\Kmat$ is of size $\numdata \numtime \times \numdata \numtime$, so precomputing it requires a total of $\numdata \numtime \sim \! 24\,500$ matvecs with $\ptomat$ and $\Gmat^*$.
Because $\Gmat^*$ has the same dimensions and block Toeplitz structure as $\ptomat^*$, matvecs with $\Gmat^*$ come at the same computational expense as matvecs with $\ptomat$ or $\ptomat^*$ ($\sim \! 0.011$~s on 24~A100~GPUs). 
Therefore, the total time to compute and store $\Kmat$ is $\numdata \numtime \cdot 0.022$~s $\sim \! 9$~m. (Note that $\bfnoise$ is assumed to be diagonal, so its computational expense in computing $\Kmat$ is negligible.) 
$\Kmat$ is a dense matrix without any obvious structure, so storing it requires storing all of its entries ($\sim \! 4.8$~GB). 
Once $\Kmat$ has been computed, it can be Cholesky-factorized ($\Kmat = \mathbf{L} \mathbf{L}^*$) which takes $\sim 1.18$~s on a single A100~GPU; the Cholesky factor $\mathbf{L}$ may overwrite $\Kmat$ and can be stored at half its memory cost ($\sim \! 2.4$~GB).\footnote{In practice, many solvers store the Cholesky factor as a lower (or upper) triangular part of the full matrix, in which case the memory cost for storing $\mathbf{L}$ equals storing $\Kmat$.}

After $\Kmat$ has been factorized, parameter uncertainties can be obtained from the posterior covariance $\bfpostcovm = \left( \Imat - \Gmat^* \Kmat^{-1} \ptomat \right) \bfpriorcovm$.
Each action of $\bfpostcovm$ requires a matvec with $\Gmat^*$ and $\ptomat$ (11~ms each), a triangular solve of the prefactorized $\Kmat$ matrix (20~ms), and a prior solve (84~ms), for a total time of $\sim \! 126$~ms.
This enables us to rapidly compute columns of $\bfpostcovm$.
However, computing the full covariance matrix is still too costly, both in memory and computation time.
Each of the $\numparam \numtime = 131\,584\,500$ columns has a storage cost of $\sim \! 1.05$~GB thus totaling $\sim \! 138$~PB for the entire $\bfpostcovm$, and computing this many matvecs would take $\numparam \numtime \cdot 126$~ms $\sim \! 192$~days.
The storage issue could be overcome by limiting ourselves to computing pointwise variances from the diagonal of $\bfpostcovm$, and the matvecs could be performed simultaneously (with each matvec performed on 24~A100 GPUs) to reduce the computation time to less than two days on a large-scale machine with 2400~A100 GPUs.
So computing parameter uncertainties may be indeed feasible for a problem of the size of this numerical example.
Unfortunately, the parameter space dimension for the CSZ digital twin is at least one order of magnitude larger ($\mc{O}(10^9)$), making a full computation of the diagonal of $\bfpostcovm$ computationally intractable.

\paragraph{Phase 3 (offline)}

As a preliminary step to computing posterior covariance of the QoIs, we precompute $\Gmat_q^* = \bfpriorcovm \ptqmat^*$.
Because $\Gmat_q^*$ is block Toeplitz, we only need to perform $\numqoi$ prior solves to compute and store  $\Gmat_q^*$.
Storage costs of $\Gmat_q^*$ are equivalent to storing $\ptqmat$ ($\sim \! 17$~GB).
The time per prior solve, performed on a single A100~GPU with Nvidia cuDSS, is $\sim \! 0.084$~s. Therefore, the total time needed to compute $\Gmat_q^*$ is $\numqoi \cdot 0.084$~s $\sim \! 1.4$~s.

To precompute the posterior covariance of the QoIs, $\bfpostcovq = \ptqmat \left( \Imat - \Gmat^* \Kmat^{-1} \ptomat \right) \Gmat_q^*$, its action needs to be evaluated $\numqoi \numtime^q = 800$ times.
Each action of $\bfpostcovq$ requires one matvec each with $\ptqmat$~(3~ms), $\Gmat^*$~(11~ms), $\ptomat$~(11~ms), and $\Gmat_q^*$~(3~ms), as well as a triangular solve of the prefactorized $\Kmat$ matrix (20~ms).
So each $\bfpostcovq$ matvec takes $\sim \! 48$~ms, leading to a total time of $\numqoi \numtime \cdot 48$~ms $\sim \! 38.4$~s to compute the full QoI covariance matrix.
For this QoI dimension of 800, the storage cost for $\bfpostcovq$ is negligible ($\sim \! 5.2$~MB).

Precomputing the d2q map $\Qmat = \ptqmat \left( \Imat - \Gmat^* \Kmat^{-1} \ptomat \right) \Gmat^* \bfinvnoise$ requires $\numqoi \numtime^q = 800$ matvecs with $\Qmat$.
Each matvec takes 56~ms, so the total time to compute $\Qmat$ is $44.8$~s, and the cost for storing the matrix is 157~MB.

\paragraph{Phase 4 (online)}

To test the inference and prediction framework, we use synthetic data generated from the PDE model. 
For the numerical example, we use a source function that describes a spatiotemporal seafloor deformation from three mixed Gaussians:
\be
   \frac{\p b(x,y,t)}{\p t} = \sum_{i=1}^{3} G^{(i)}(x,y,t) ,
   \label{eq:source-example}
\ee
where
\be
   G^{(i)}(x,y,t) = \left\{ 
   \begin{aligned}
   & A^{(i)} \exp(-\left(\frac{x-X_{\text{c}}^{(i)}}{X_{\text{r}}^{(i)}}\right)^2-\left(\frac{y-Y_{\text{c}}^{(i)}}{Y_{\text{r}}^{(i)}}\right)^2) \frac{\pi}{2 T_{\text{r}}^{(i)}} \sin(\frac{\pi t}{T_{\text{r}}^{(i)}}) & \! \! \! \! \!, \quad t \le T_{\text{r}}^{(i)}, \\
   & 0 & \! \! \! \! \!, \quad t > T_{\text{r}}^{(i)} \, .
   \end{aligned} \right.
   \label{eq:source-example-Gaussians}
\ee
The parameters for these Gaussians are given in Table~\tab{seafloor-gaussians}.
While each of the Gaussians varies smoothly in both space and time, the second and third ($i=2,3$) introduce higher-frequency components (both spatially and temporally) into the solution than the seafloor motion due to the first Gaussian ($i=1$).

\begin{table}[htb]
   \centering
   \caption{Parameters of the Gaussians defining synthetic seafloor deformation for the numerical example.}
   \label{tab:seafloor-gaussians}
   \begin{tabular}{llllll}
      \toprule
      Symbol & Description & & Value & & Unit \\
      \cmidrule{3-5}
      & & $i=1$ & $i=2$ & $i=3$ & \\
      \midrule
      $A^{(i)}$ & Amplitude & 4 & 1 & -0.5 & m \\
      $T_{\text{r}}^{(i)}$ & Rise time & 20 & 10 & 10 & s \\
      $X_{\text{r}}^{(i)}$ & Rise width in $x$ & 16 & 4 & 4 & km \\
      $Y_{\text{r}}^{(i)}$ & Rise width in $y$ & 32 & 4 & 8 & km \\
      $X_{\text{c}}^{(i)}$ & Center in $x$ & 64 & 64 & 70 & km \\
      $Y_{\text{c}}^{(i)}$ & Center in $y$ & 64 & 88 & 56 & km \\
      \bottomrule
   \end{tabular}
\end{table}

Noisy data are generated from the true data (obtained from the PDE model) by adding Gaussian noise. In particular, for each sensor, the magnitude of the added Gaussian noise is taken as a certain percentage (from now on referred to as ``noise level'') of the infinity norm (over all time steps) of the sensor's true data.
Given these noisy observations, the total online compute time to solve the inference problem for $131\,584\,500$ parameters (via \eq{real-time-parameter}) is $\sim 53$~ms on 24 A100 GPUs.\footnote{Neglecting cost for I/O, the total compute time is composed of $\sim 11$~ms per $\Gmat^*$ and $\ptomat$, and $\sim 20$~ms for solving prefactorized $\Kmat$.}
Solving the prediction problem for $800$ QoIs can then be done either through push-forward of the inferred parameters via the p2q map ($\sim 3$~ms on 24 A100 GPUs) or, more efficiently, directly from the noisy observations (via \eq{real-time-qoi}) using the precomputed d2q map ($\sim 1$~ms on a laptop).
In other words, the full end-to-end digital twin, from observations via inference to prediction, with the full-order model and including uncertainties, is realized within fractions of a second for a test configuration of nearly 132~M parameters.

\begin{table}[H]
   \centering
   \caption{Compute times and storage requirements for offline computation (Phases~1--3) and online computation (Phase~4) of the Bayesian inference and prediction framework for a numerical example with 132~M parameters.}
   \label{tab:time-and-memory}
   \begin{tabular}{lllll}
      \toprule
      Phase & Task & Time & Memory & System \\
      \midrule
      1 & Precompute p2o map $\ptomat$ & 6.5~h & 52~GB & 3584 Sapphire Rapids Cores \\
         & (49 adjoint PDE solves) \\
         & Precompute p2q map $\ptqmat$ & 2.1~h & 17~GB & \\
         & (16 adjoint PDE solves) \\
      \cmidrule{1-5}
      2 & Form $\Gmat^*$ (49 prior solves) & 4.2~s & 52~GB & 24 Nvidia A100 (40~GB) GPUs \\
         & Form $\Gmat_q^*$ (16 prior solves) & 1.4~s & 17~GB & \\
         & Compute and factorize $\Kmat$ & 10~m & 4.8~GB & \\
      \cmidrule{1-5}
      3 & Compute QoI cov.\ $\bfpostcovq$ & 1~m & 5.2~MB & 24 Nvidia A100 (40~GB) GPUs \\
         & Precompute d2q map $\Qmat$  & 1~m & 0.16~GB & \\
      \cmidrule{1-5}
      4 & Inference of $\bfmmap$ & 57~ms & 1.1~GB & 24 Nvidia A100 (40~GB) GPUs \\
         & Prediction of $\bfqmap$ & 1~ms & 0.1~MB & Laptop \\
      \bottomrule
   \end{tabular}
\end{table}

Table~\tab{time-and-memory} summarizes the computational costs for each phase of the framework, both in storage requirements and compute times (obtained on systems as specified in the table), for the test configuration.
While the spatial domain of this numerical example was about $20 \times$ smaller than the dimensions of the full CSZ digital twin, the computational resources can be scaled up accordingly.
Importantly, all of the large-scale computations done offline---particularly the PDE solves in Phase~1---scale effectively to larger compute clusters.

\begin{table}[htb]
   \centering
   \caption{Relative errors of the inferred parameters $\bfmmap$, the predicted QoIs $\bfqmap$, and the pressures at sensor locations reconstructed from $\bfmmap$ via the p2o map $\ptomat$, obtained via the real-time inference and prediction framework from synthetic observations with varying levels of relative additive noise.}
   \begin{tabular}{llll}
      \toprule
      Noise level & Inferred parameters & Predicted QoIs & Pressure reconstruction at sensors \\
      \midrule
      2\% & 0.0776 & 0.0108 & 0.0195 \\
      4\% & 0.0836 & 0.0167 & 0.0397 \\
      6\% & 0.0948 & 0.0161 & 0.0596 \\
      \bottomrule
   \end{tabular}
   \label{tab:relative-error}
\end{table}

For the particular source function given in \eq{source-example}--\eq{source-example-Gaussians} and three examples of synthetic observations with varying additive noise levels (2\%, 4\%, 6\%), Table~\tab{relative-error} shows the relative errors observed for the mean of the inferred parameters, the mean of the predicted QoIs, and the pressures at sensor locations reconstructed from the mean of the inferred parameter field via the p2o map.
The relative error of the reconstructed pressures, $\| \ptomat \bfmmap - \ptomat \bfmtrue \| / \| \ptomat \bfmtrue \|$, corresponds roughly to the noise level; this is expected because the highly underdetermined inverse problem ($\numparam \gg \numdata$) closely fits the sparse and noisy data, even with the smoothing prior.
The relative errors of the inferred mean of the parameter field, $\| \bfmmap - \bfmtrue \| / \| \bfmtrue \|$ ($\sim$~7--10\%), are larger than those of the predicted mean of the QoIs, $\| \bfqmap - \bfqtrue \| / \| \bfqtrue \|$ ($\sim$~1--2\%).
Since the QoIs are predicted via the inference of the parameters, this result suggests that some components of the error of the inferred parameter field are not forward propagated to the QoIs via the p2q map.

\begin{figure}[htb]
   \centering
   \begin{subfigure}[t]{0.32\textwidth}
      \includegraphics[width=\textwidth,clip,trim={0 0 0 0}]
      {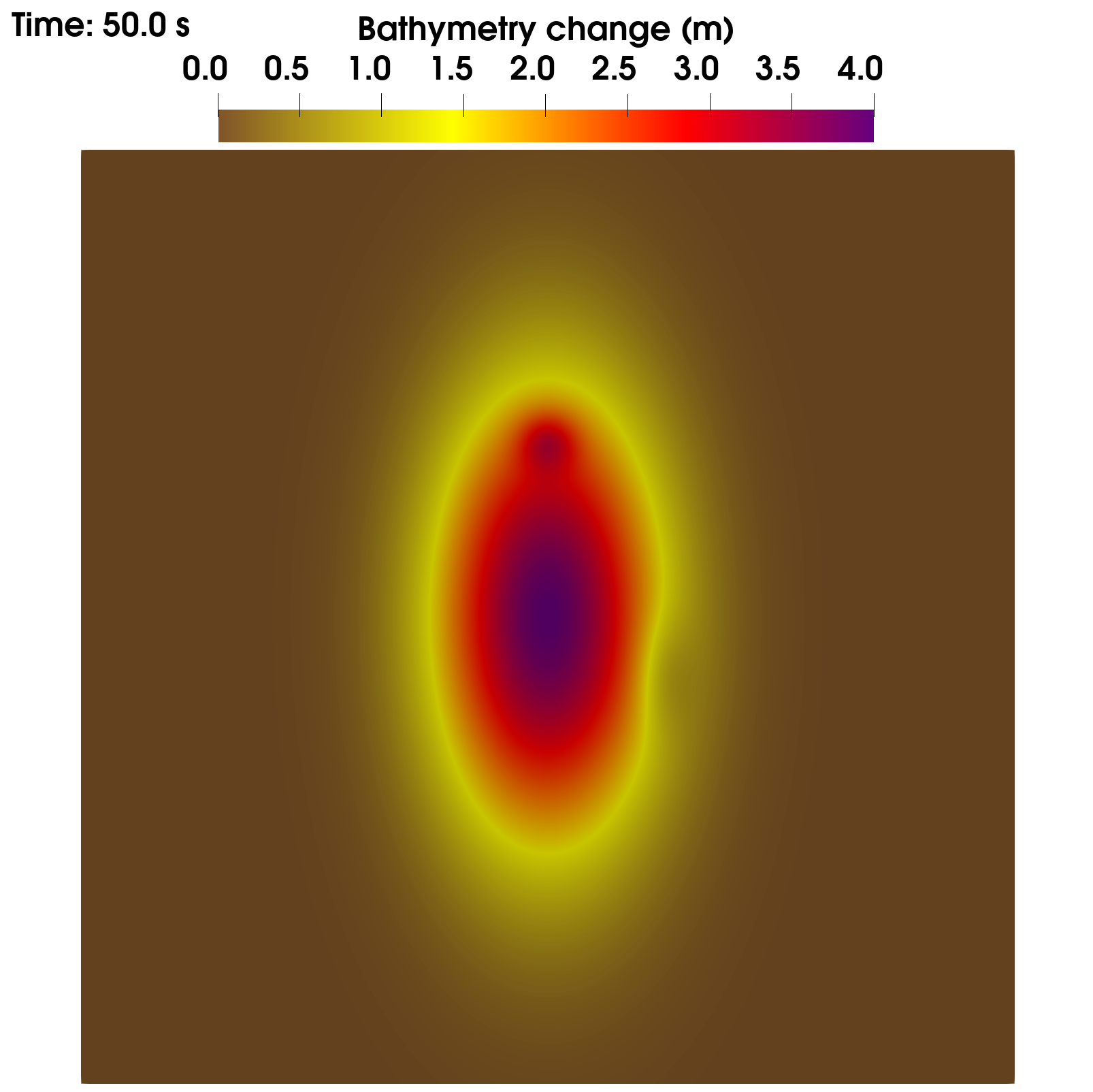}
   \end{subfigure}
   \begin{subfigure}[t]{0.32\textwidth}
      \includegraphics[width=\textwidth,clip,trim={0 0 0 0}]
      {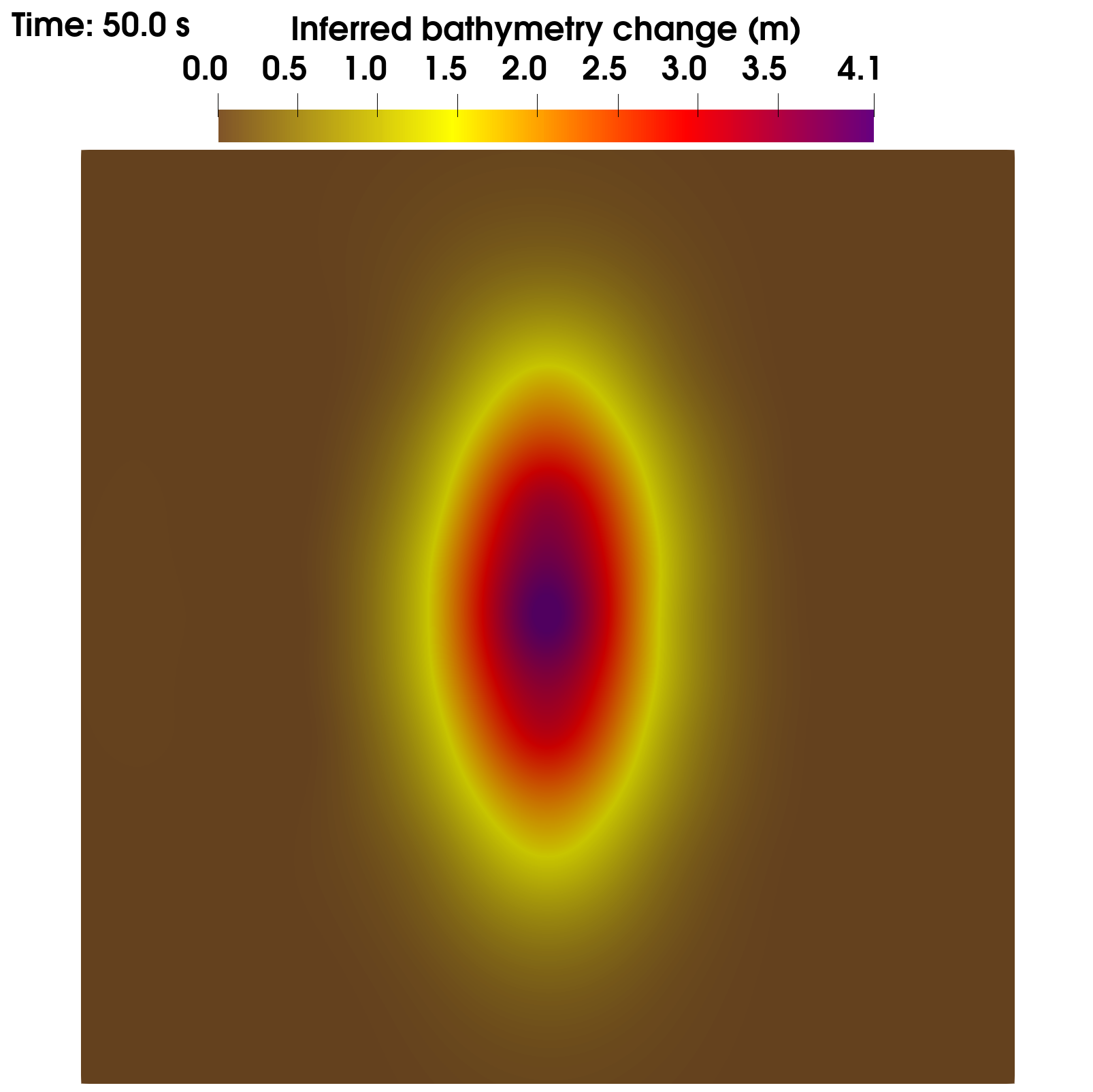}
   \end{subfigure}
   \begin{subfigure}[t]{0.32\textwidth}
      \includegraphics[width=\textwidth,clip,trim={0 0 0 0}]
      {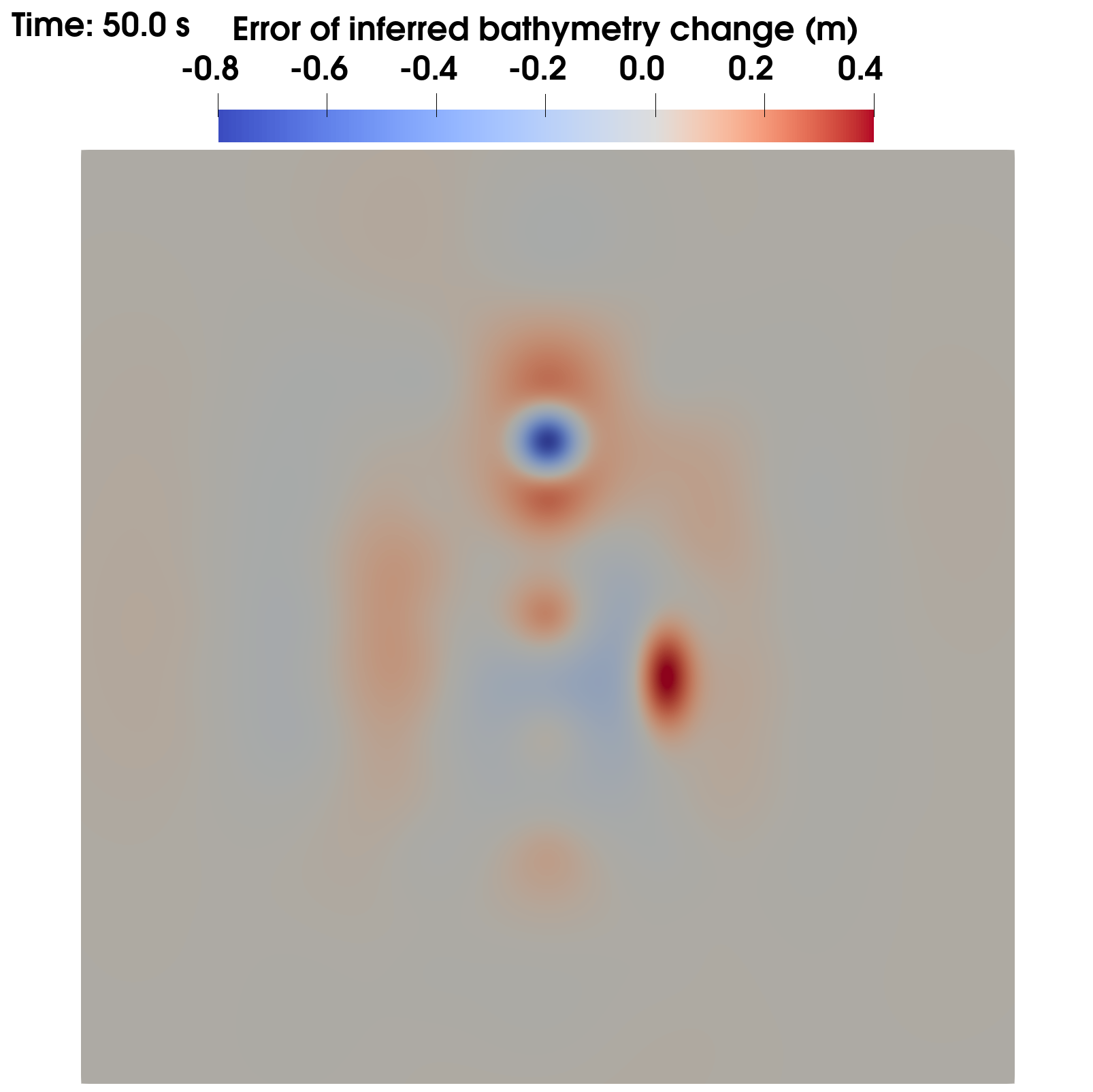}
   \end{subfigure}
    \caption{Total seafloor displacement (time-integrated parameter field). Left: True field; Middle: Inferred mean from synthetic observations with 6\% relative additive noise; Right: Error between true and inferred fields.}
   \label{fig:example-parameter}
\end{figure}

\begin{figure}[htb]
   \centering
   \begin{subfigure}[t]{0.32\textwidth}
      \includegraphics[width=\textwidth,clip,trim={0 0 0 0}]
      {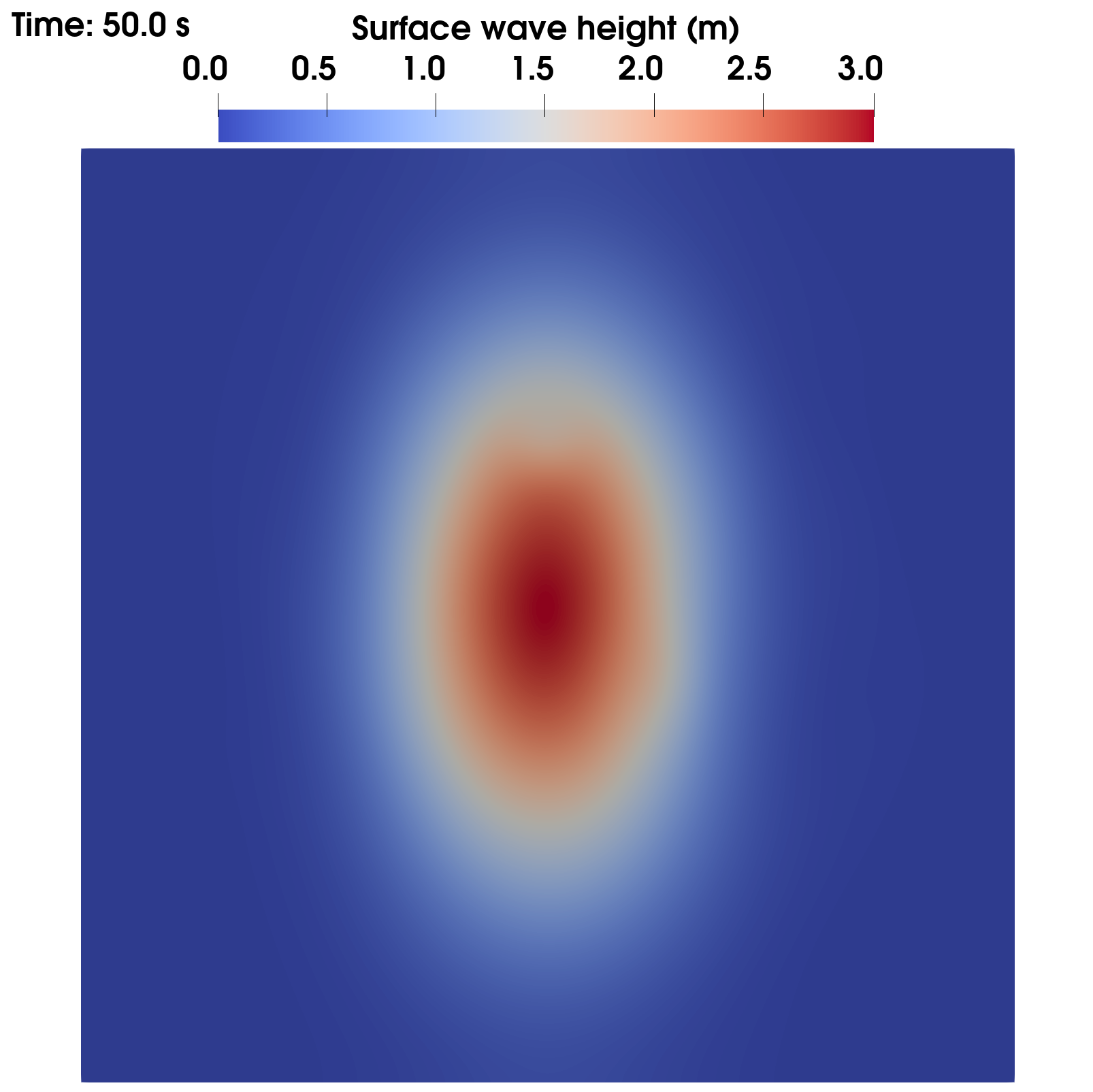}
   \end{subfigure}
   \begin{subfigure}[t]{0.32\textwidth}
      \includegraphics[width=\textwidth,clip,trim={0 0 0 0}]
      {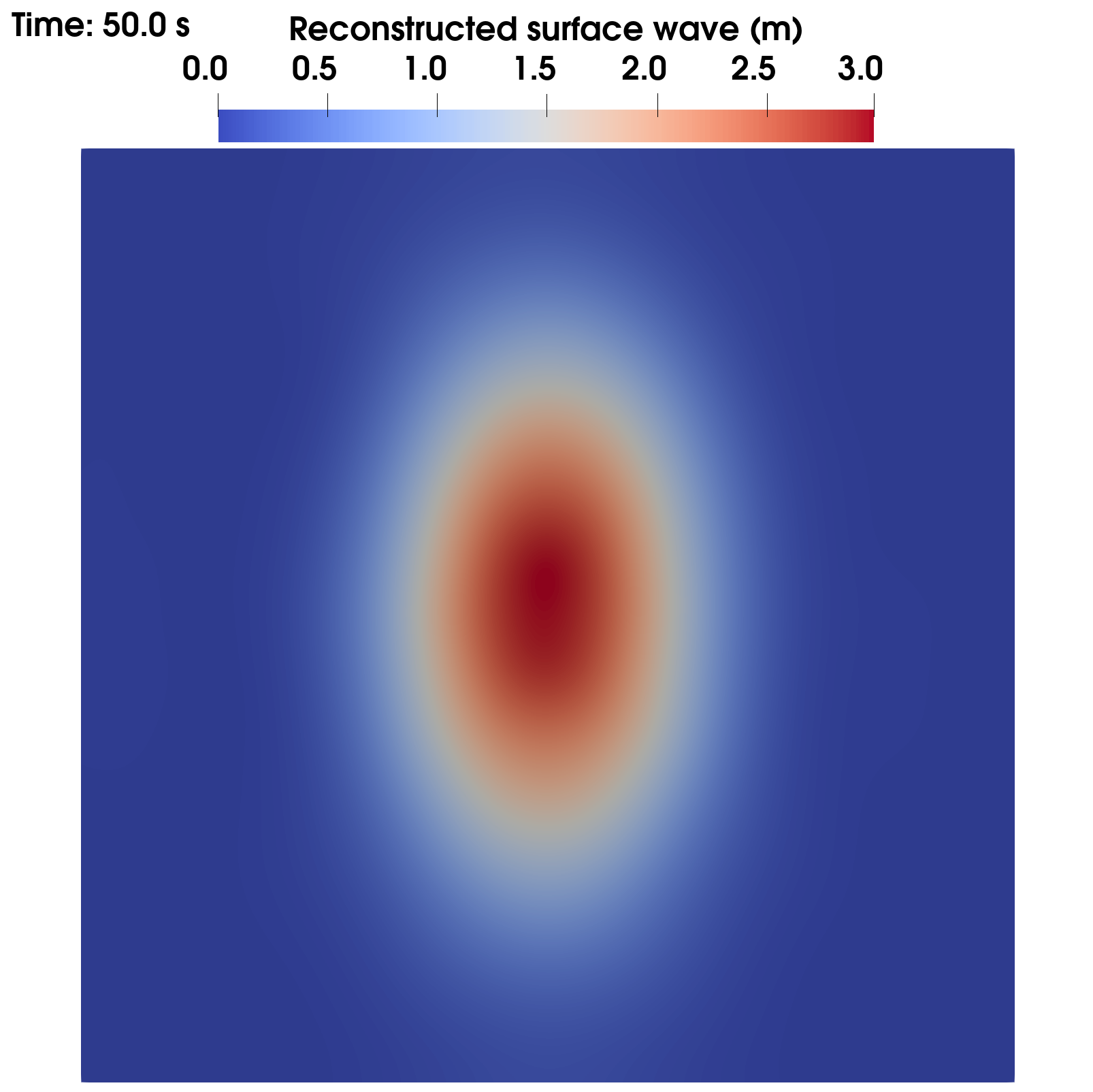}
   \end{subfigure}
   \begin{subfigure}[t]{0.32\textwidth}
      \includegraphics[width=\textwidth,clip,trim={0 0 0 0}]
      {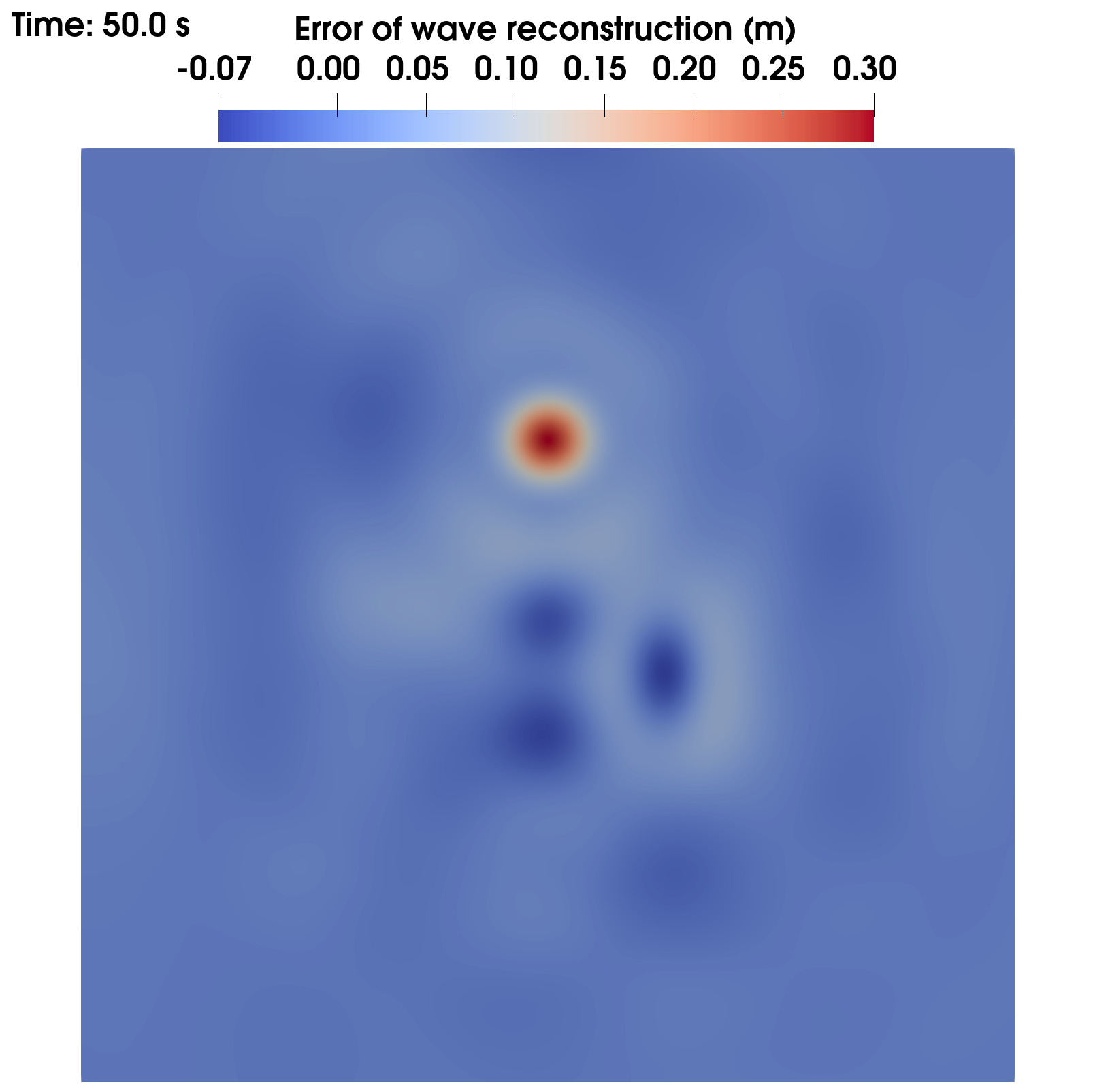}
   \end{subfigure}
   \caption{Surface gravity wave. Left: True field; Middle: Reconstructed field based on inferred mean of seafloor motion (from synthetic observations with 6\% relative additive noise); Right: Error between true and reconstructed surface gravity wave based on inferred mean of seafloor motion.}
   \label{fig:example-surface-wave}
\end{figure}

\begin{figure}[htb]
   \centering
   \begin{subfigure}[t]{0.49\textwidth}
      \includegraphics[width=\textwidth,trim={0 0 0 0}]
      {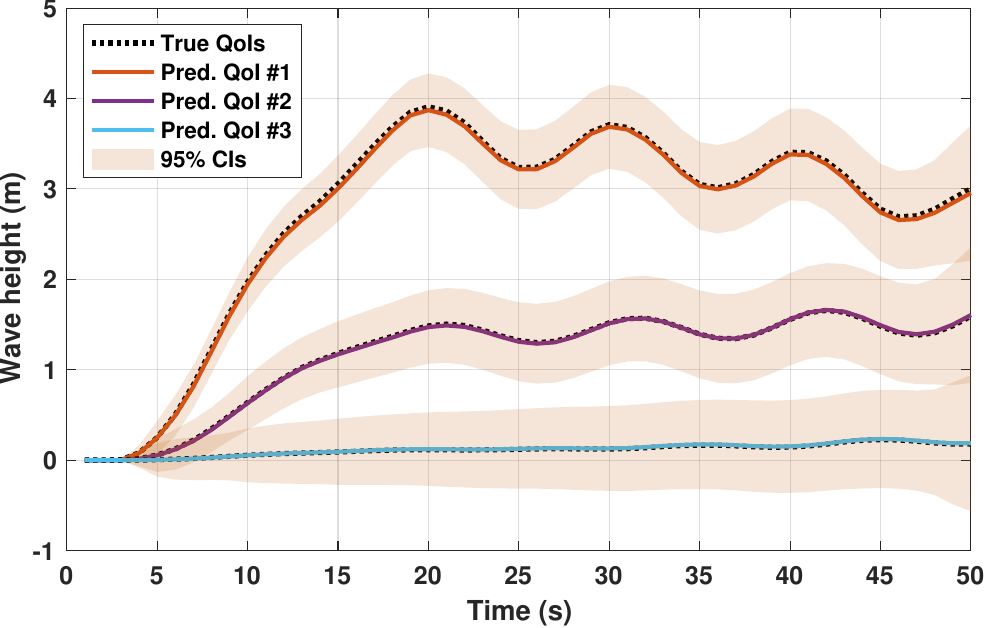}
   \end{subfigure}
   \begin{subfigure}[t]{0.49\textwidth}
      \includegraphics[width=\textwidth,trim={0 0 0 0}]
      {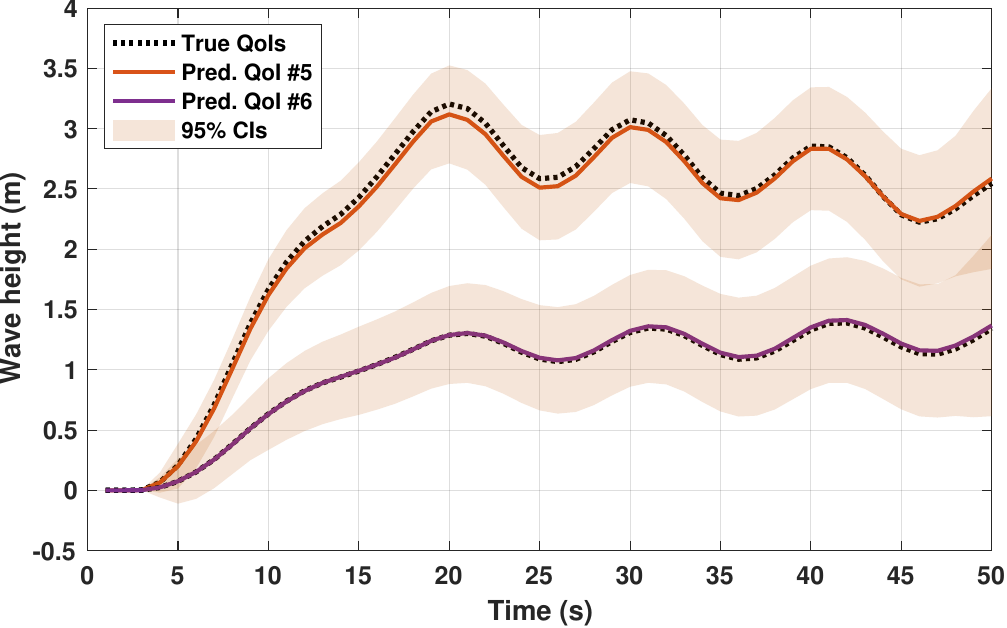}
   \end{subfigure}
   \caption{True and predicted QoIs (tsunami wave heights at specified locations) with uncertainties depicted in the form of 95\% credible intervals (CIs). Left: QoI locations \#1--3; Right: QoI locations \#5--6 (cf.\ Figure~\fig{domain}).}
   \label{fig:example-qoi}
\end{figure}

Figure~\fig{example-parameter} illustrates that the components of the parameter field that are not inferred correctly predominantly correspond to the higher-frequency spatial components.
These rough components cannot be reconstructed because they are not well informed by the sparse data, a fundamental limitation of an ill-posed inverse problem~\cite{ghattas2021learning}.
On the other hand, the lower-frequency modes of the seafloor motion are inferred with good accuracy.
A similar observation can be made for the surface gravity wave reconstructed based on the inferred mean of the seafloor motion, which is depicted in Figure~\fig{example-surface-wave}.
However, the rough components of the relative error for the surface gravity wave are smaller than those for the parameter field, suggesting that the corresponding higher-frequency modes are dampened by the p2q map.
In other words, the tsunami formation appears to be driven mostly by the lower-frequency components of the seafloor motion, making relatively accurate QoI predictions possible even when the inferred parameters have larger relative errors (cf.\ Table~\tab{relative-error}) and those errors are dominated by the rough components.
The inferred mean of some of the predicted QoIs (cf.\ Figure~\fig{domain}), as well as their respective 95\% credible intervals, are shown in Figure~\fig{example-qoi}.

Finally, we compare our method to matrix-free conjugate gradients (CG), a conventional state-of-the-art method for solving large-scale inverse problems. In this numerical example, the effective rank of the Hessian (Figure~\ref{fig:eigenvalues}) is nearly of the order of the data dimension (here $24\,500$), and thus we can expect CG to take $\mathcal{O}(24\,500)$ iterations. Since each iteration requires the action of $\ptomat$ and $\ptomat^*$ on a vector, i.e.~a pair of forward/adjoint wave propagations ($\sim 16$~minutes on $3\,584$ Sapphire Rapids cores), we would expect CG to take $\mathcal{O}(272)$ days for this model problem. And this is for computing $\bfmmap$ via \eq{mmap-problem}; evaluating the uncertainty of the inverse solution (also needed for the posterior of the QoI) via the inverse of the Hessian \eq{Hessian} is entirely intractable. In comparison, using our method the online inference and prediction stage takes $58$~ms (on $24$ A100 GPUs), a 400-million-fold speedup. Factoring in the one-time offline cost of 65 adjoint wave propagations, the method results in a $\sim 750 \times$ reduction in PDE solutions if the inverse problem is solved just once.

\section{Conclusions and Ongoing Work}
\label{sec:conclusions}

By exploiting the inherent autonomous structure of the governing
forward PDE model, we have developed a digital twin framework for
real-time Bayesian inversion at extreme scale while remaining faithful
to high-fidelity, physics-based models.  As a core component of this
framework, we employ an FFT-based algorithm to compute Hessian matvecs
at speeds over four orders of magnitude faster than those of
conventional methods and circumvent the need for PDE solution at the
time of inference. We re-emphasize that these FFT-based Hessian matvecs are 
exact and do not incorporate any approximations. In addition to making real-time 
Hessian matvecs possible, this algorithm can also be used to accelerate offline
portions of the computation, which includes precomputing an efficient
representation of the inverse operator in the data space.  The digital
twin framework is extended to the goal-oriented setting to enable
real-time predictions of QoIs under uncertainty.

As a target application, we utilized this framework for a
proof-of-concept tsunami early warning system that infers the
earthquake-induced seafloor motion from near-field observations of the
acoustic pressure on the seafloor and an acoustic--gravity PDE model,
followed by a forward prediction of tsunami wave heights at specified
locations with associated uncertainties.  The entire end-to-end
inference, prediction, and uncertainty quantification for a system
with over one hundred million parameters was computed exactly in
fractions of a second.  For this tsunami inversion, the results
indicate that sparse observations of the pressure transients in the
near field are sufficient to accurately infer the smooth components of
the seafloor motion and that these smooth components appear to account
well for the surface gravity wave formation.

In ongoing work, we are leveraging this framework to construct a
digital twin for tsunami forecasting from potential rupture events in
the CSZ.  Implementation of the CSZ digital twin requires offline
computations of extreme-scale adjoint wave propagation PDE solutions
and online inference of more than one billion parameters (i.e.~at
least one order of magnitude larger than the model problem solved
in this paper) \cite{henneking2025bell}.  Beyond tsunami early
warning, the core components of our digital twin framework are
applicable to many other classes of dynamical system inverse problems,
including nuclear testing treaty verification, contaminant transport,
and atmospheric CO2 monitoring.  Extending the framework to construct
digital twins for these problems, and others, is ongoing work.

\section*{Acknowledgements}
This research was supported in part by DOE ASCR grant DE-SC0023171 and DOD MURI grants FA9550-21-1-0084 and FA9550-24-1-0327.
Supercomputing resources were provided by the Texas Advanced Computing Center (TACC) at UT Austin on its \emph{Frontera}, \emph{Lonestar6}, \emph{Stampede3}, and \emph{Vista} systems, as well as by NERSC award ALCC-ERCAP0030671 at NERSC on its \emph{Perlmutter} system.
This material is based upon work supported by the National Science Foundation Graduate Research Fellowship under Grant No.~DGE 2137420.
We also thank Dr.~Milinda Fernando from The University of Texas at Austin for insightful discussions.


\printbibliography[heading=bibintoc]


\end{document}